\documentclass[lettersize,journal]{IEEEtran}

\usepackage{amsmath,amssymb,amsfonts}

\usepackage{algorithmic}
\usepackage{algorithm}
\usepackage{array}
\usepackage[caption=false,font=normalsize,labelfont=sf,textfont=sf]{subfig}
\usepackage{textcomp}
\usepackage{stfloats}
\usepackage{url}
\usepackage{verbatim}
\usepackage{graphicx}
\usepackage{cite}
\DeclareMathOperator{\E}{E}
\DeclareMathOperator{\diag}{diag}
\DeclareMathOperator{\var}{Var}

\newtheorem{remark}{Remark}
\newtheorem{theorem}{Theorem}
\newtheorem{defi}{Definition}
\newtheorem{proposition}{Proposition}
\newtheorem{property}{Property}

\begin{document}

\title{Multi-sensor Distributed Fusion Estimation for $\mathbb{T}_k$-proper Factorizable  Signals  in Sensor Networks with Fading Measurements}

\author{Rosa M.  Fern\'andez-Alcal\'a, Jos\'e D. Jim\'enez-L\'opez, Jes\'us Navarro-Moreno, and Juan C.  Ruiz-Molina
\thanks{The authors are with the Department of Statistics and Operations Research, University of Ja\'en, 23071 Ja\'en, Spain. E-mails: {rmfernan,jdomingo,jnavarro,jcruiz}@ujaen.es. Tel:+34 953212449, Fax:+34 953212034.}
\thanks{This work is part of the I+D+i  project PID2021-124486NB-I00, funded by MICIU/AEI/10.13039/501100011033/ and ERDF/EU. }
}

\markboth{IEEE Transactions  on Signal and Information Processing over Networks,~Vol.~XX, No.~Y,  MONTH 2025}%
{Fern\'andez-Alcal\'a \MakeLowercase{\textit{et al.}}: Multi-sensor Distributed Fusion Estimation for $\mathbb{T}_k$-proper Factorizable  Signals  in Sensor Networks with Fading Measurements}



\maketitle

\begin{abstract}
The challenge of distributed fusion estimation is investigated for a  class of four-dimensional (4D) commutative hypercomplex signals that  are $\mathbb{T}_k$-proper factorizable,  within the framework of multiple-sensor networks with  different fading measurement rates. 
The fading effects affecting each sensor's measurements are modeled as a stochastic variables with known second-order statistical properties. The estimation process is conducted exclusively based on these second-order statistics.  Then, by exploiting the $\mathbb{T}_k$-properness property within a  tessarine framework,  the dimensionality of the problem is significantly reduced. This reduction in dimensionality   enables the development of distributed fusion filtering,  prediction, and smoothing algorithms that entail lower computational effort compared with real-valued approaches.
  The performance of the suggested algorithms is assessed through numerical experiments under various uncertainty conditions and $T_k$-proper contexts.   Furthermore,  simulation results confirm  that $\mathbb{T}_k$-proper  estimators outperform their quaternion-domain counterparts,  
  underscoring their practical advantages. These findings highlight the potential of 
$\mathbb{T}_k$-proper estimation techniques for improving multi-sensor data fusion in applications where efficient signal processing is essential.
\end{abstract}


\section{Introduction}
\label{sec:introduction}
\IEEEPARstart{M}{ulti-sensor}  data fusion involves combining information from various sources to produce a more accurate, reliable and comprehensive understanding of a situation.  This theory emerged in the 1970s to meet military needs. Nevertheless,  advances in sensor technology have made  it more accessible,  and nowadays,  the use of multiple sensor systems 
has spread to a wide variety of application fields, including healthcare \cite{Singh2024},   robotics \cite{Kurkin2017},  navigation and location \cite{Zewge2023, Yan2019}, video and image processing \cite{Huang2020},  and  communication networks \cite{ Liping2021}.  

In this context, extensive literature exists on fusion estimation in multi-sensor systems,   in which various fusion strategies have been applied (see, e.g.,  \cite{Liping2021,  
Hu2020},  and references therein).   In broad terms,  centralized and distributed fusion techniques are  the most commonly used fusion structures.  In the former,  measurement information from  individual sensors is transmitted  to a fusion center for estimation.  This approach typically yields  optimal estimators when all the sensors operate without any faults,  but they can suffer from a high computational cost.
In contrast, distributed fusion employs decentralized methods where measurement data are processed independently  at each sensor, and the resulting local estimators are then fused to produce a more accurate overall estimate.  Compared to the centralized fusion architecture,  this methodology enhances scalability and  robustness,   thereby minimizing the risk of a single point of failure,  and also it reduces overall computational load. 
These advantages have generated significant interest in research on distributed fusion techniques (see, e.g.,  \cite{
Hu2020,Li2022}).

Moreover,  in practical applications,  network systems often operate in unreliable environments where stochastic uncertainties 
frequently  arise, potentially  degrading system performance .  A significant source of uncertainty in such systems is fading measurements. They reflect phenomena characterized by unpredictable fluctuations  in signal strength caused by  aging sensors, external interferences,   or faulty communication channels \cite{
Jin2022}.   

Most of  the  distributed fusion estimation algorithms in multi-sensor systems with fading measurements rely on state-space models \cite{Jin2022, 
Sun2018, 
Hu2024,Jin2024}.  
Among the different fusion architectures,  Kalman consensus filtering has been used in \cite{Jin2022,  Jin2024} for linear and nonlinear systems with multiple uncertainties,    
matrix-weighted fusion estimation appears in \cite{Sun2018},  whereas
 covariance intersection fusion methods have been implemented  in \cite{
 Hu2024}. 
 Alternatively,   the use of second-order statistical information has proven to be a key strategy for addressing fusion estimation problems.  In particular,  by assuming factorizable covariances,  matrix-weighted-based distributed fusion linear estimators have been provided in \cite{GarciaLigero2020}  for systems characterized by mixed uncertainties 
  This structural assumption is very general and handles both stationary and non-stationary signals,  offering  a comprehensive  framework for addressing a wide range of real-world situations
 \cite{ Caballero2020b}. 
Moreover,   unlike state-space models, which require a detailed characterization of system dynamics and measurement processes, the factorizable covariances approach directly  targets the covariance structure of the processes involved. This makes it more flexible and easier to use,  particularly in scenarios where detailed system models are unavailable or impractical.

In addition to traditional methods, the use of four-dimensional (4D) hypercomplex algebras has attracted attention for solving multi-sensor fusion estimation problems (see,  e.g.,  \cite{Talebi2016, Talebi2020, Jimenez2021, Fernandez2023, Farkas2024}). These algebras  extend the concept of complex numbers by adding three imaginary units to the real component, and  offer a versatile framework for modeling three- and four- dimensional phenomena  \cite{Ortolani2019,  Valle2020, Valle2021,  Cariow2021,  Grassucci2021,  Vieira2022,  Guizzo2023,  
Clive2024, Cariow2024,  Kosal2024a,  Kosal2024b}
Quaternions are among the most prominent hypercomplex structures in signal processing,  demonstrating their effectiveness due to their distinctive ability to represent rotations and handle multi-dimensional data efficiently \cite{Guizzo2023, 
Clive2024}. 
Even though  quaternions are the only 4D hypercomplex algebra without zero divisors,  their non-commutative nature poses challenges in certain applications.   More recently,  tessarines have gained attention  as a type of commutative 4D algebra that offers  significant advantages over quaternions: Their simpler structure allows extending 
results from real and complex fields to four-dimensional space \cite{Ortolani2019},   they are computationally more efficient  in many scenarios \cite{Pei2004, Pei2008},  and  have also demonstrated superior performance in various domains,  including image processing,  signal processing and pattern recognition  \cite{Melegy2022, Guo2024}.  In practical applications, tessarines have been applied in many diverse areas, including electromagnetic theory,  color images processing,  images classification,  and neural networks  \cite{Melegy2022, Guo2024, Gai2024}.

Furthermore,  when working with hypercomplex signals, it is possible to define certain properness characteristics related to the second-order statistics that  lead to  dimensionality reduction and, consequently, computational simplifications in the estimation algorithms devised -- advantages not attainable within a real-valued framework.  Specifically,   while the optimal method for processing 4D hypercomplex signals is the widely linear (WL) processing,  that means operating on a 4D vector, under properness conditions, the same optimal solutions can be achieved by operating on a vector of  half or  quarter the original dimension \cite{Nitta2019, Navarro2021}. 
This fact makes hypercomplex signals particularly useful for applications in multi-sensor fusion, where the ability to obtain estimates of the signals with a lower computational cost  is key. 
So far, in the tessarine field,  by assuming  $\mathbb{T}_k$-properness conditions ($k=1,2$),   
 low-cost centralized and distributed fusion estimation algorithms  have been devised for  multi-sensor state-space systems 
 with different uncertainties (missing measurements,  random delays,  and packet dropouts)    \cite{Jimenez2021, Fernandez2023}.  Additionally,  the concept of factorizable kernels has been extended to the tessarine domain,  leading to the  introduction of widely factorizable signals. These signals are characterized by the property that the correlation matrix of the 4D augmented vector —composed of the signal and three auxiliary functions—forms a factorizable kernel. This class of signals is quite general and can represent both stationary and nonstationary signals \cite{Fernandez2019}.

Therefore,  this paper investigates multi-sensor systems with fading measurements on the class of tessarine signals which are both widely factorizable  and $\mathbb{T}_k$-proper.  
In this context,  the  filtering, prediction, and smoothing estimation problems are addressed.
Then, taking the computational advantages of the $\mathbb{T}_k$-proper signals and using solely
second-order statistical information, distributed fusion estimation algorithms are devised for all three  estimator types. These algorithms exhibit two main advantages: 1) applicability in scenarios where  a state-space model is not directly accessible, and 2) reduced computational complexity compared to their real-valued counterparts.

The structure of the paper is as follows.  Section 2 establishes the notation used throughout and  provides a review of  the fundamentals of the  tessarine algebra.  Section 3 formulates the estimation problem and explores the impact of $\mathbb{T}_k$-properness on multi-sensor systems with fading measurements.  In Section 4,  the distributed fusion linear estimation problem is tackled.  Firstly,  based on the available second-order statistical information,  linear minimum mean squared error (MMSE) algorithms are designed at each sensor for the filtering, prediction,  and smoothing estimators. Then, by combining these local estimators through a weighted linear approach, distributed fusion estimation algorithms are proposed to compute the overall estimators and their associated mean squared errors. Section 5 includes numerical experiments  that illustrate the effectiveness and practical relevance of the methods proposed.  Finally,  Section 6 presents the concluding remarks.

\section{Basic Knowledge}

This section establishes the notational conventions adopted throughout the paper and provides a brief overview of the key concepts and properties of tessarine algebra.

\subsection{Notation}

In general, the symbols and notation used in this paper follows conventional standards.  Specifically,  scalar values are represented by  regular lowercase letters  vectors by lowercase characters,  and matrices by bold uppercase letters.  In addition,  bold upper case italicized letters stand for some specific matrices.  
Furthermore,  the symbols listed below have been employed: 

\begin{tabular}{ll}

$\boldsymbol{I}_{n}$ & {$n\times n$-dimensional identity matrix}\\
$\mathbf{0}_{n\times q}$ & {$n\times q$-dimensional matrix of zeros}\\
$\boldsymbol{1}_n$ & {$n$-dimensional vector of ones}\\
$\boldsymbol{0}_n$ & {$n$-dimensional vector of zeros}\\
 
 $(\cdot)^*$ & {Tessarine conjugate}\\
$(\cdot)^{\mathtt{T}}$ & {Transpose}\\
 $(\cdot)^{\mathtt{H}}$  & {Hermitian transpose}\\

$(\cdot)_{\rm{r}}$ & {Real part of a tessarine}\\
$(\cdot)_{{\nu}}$ & {Imaginary part ($\nu=\rm{i},\rm{j},\rm{k}$) of a tessarine}\\

$\mathbb{Z}$ & Integer field\\
$\mathbb{R}$ & Real field\\
$\mathbb{T}$ & Tessarine field\\

$\mathbf{A}\in \mathbb{R}^{n\times q}$ & 
 $n\times q$-dimensional  real matrix,\\
$\mathbf{r}\in
\mathbb{R}^{n}$  & $n$-dimensional real vector\\
$\mathbf{A}\in \mathbb{T}^{n\times q}$ &   $n\times q$-dimensional tessarine matrix\\
$\mathbf{r}\in
\mathbb{T}^{n}$  & $n$-dimensional tessarine vector\\

$E[\cdot]$ & Expectation operator \\

$\diag(\cdot)$& Diagonal (or block diagonal) matrix with \\
&  specified diagonal elements \\

$\delta_{t,s}$ &  Kronecker delta function\\

$\circ$ & Hadamard product\\
$\otimes$ &  Kronecker product\\

\end{tabular}

\medskip

 In the context of this paper,  unless explicitly mentioned otherwise,  all random variables are taken to be of zero mean.

\subsection{Tessarine algebra}

An $n$-dimensional ($n$D) tessarine stochastic signal $\mathbf{x}(t)\in \mathbb{T}^{n}$ is defined as 
\begin{equation*}\label{q1}
\mathbf{x}(t)=\mathbf{x}_{\mathrm{r}}(t)+\imath \mathbf{x}_{\imath}(t)+\jmath\mathbf{x}_{\jmath}(t)+\kappa\mathbf{x}_{\kappa}(t),\qquad t\in \mathbb{Z},
\end{equation*}
 where $\mathbf{x}_{\mathrm{r}}(t)$, $\mathbf{x}_{\imath}(t)$, $\mathbf{x}_{\jmath}(t)$ and $\mathbf{x}_{\kappa}(t)$ are $n$D real-valued random signals,  and $\{\imath,\jmath,\kappa \}$ are imaginary tessarine units with the multiplication rules: $\jmath^2=1$,  $
 \imath^2=\kappa^2=-1$,  $\imath \jmath=\kappa$,  $ 
\jmath \kappa =\imath$,  $\kappa \imath =- \jmath$.

 The conjugate of $\mathbf{x}(t)\in \mathbb{T}^n$ is given by
$\mathbf{x}^{*}(t)=\mathbf{x}_{\mathrm{r}}(t)-\imath \mathbf{x}_{\imath}(t)+\jmath \mathbf{x}_{\jmath}(t)-\kappa\mathbf{x}_{\kappa}(t)$.  Additionally,   two auxiliary tessarine vectors are considered:
 $\mathbf{x}^{\imath}(t)=\mathbf{x}_{\mathrm{r}}(t)+\imath \mathbf{x}_{\imath}(t)-\jmath \mathbf{x}_{\jmath}(t)-\kappa\mathbf{x}_{\kappa}(t)$, and
$\mathbf{x}^{\kappa}(t)=\mathbf{x}_{\mathrm{r}}(t)-\imath \mathbf{x}_{\imath}(t)-\jmath \mathbf{x}_{\jmath}(t)+\kappa\mathbf{x}_{\kappa}(t)$.

To entirely encompass the complete second-order statistical information of $\mathbf{x}(t)$,   the $4n$D  augmented signal vector $\bar{\mathbf{x}}(t)=[\mathbf{x}^{\mathtt{T}}(t), \mathbf{x}^{\mathtt{H}}(t), \mathbf{x}^{\imath^{\mathtt{T}}}(t),\mathbf{x}^{\kappa^{\mathtt{T}}}(t)]^{\mathtt{T}}$ must be considered.   Moreover,  $\bar{\mathbf{x}}(t)$ satisfies the following relationship with the corresponding real-valued vector $\mathbf{x}^r(t)=[\mathbf{x}_{\mathrm{r}}^{\mathtt{T}}(t),\mathbf{x}_{\imath}^{\mathtt{T}}(t),\mathbf{x}_{\jmath}^{\mathtt{T}}(t),\mathbf{x}_{\kappa}^{\mathtt{T}}(t)]^\mathtt{T}$:
\begin{equation*}\label{q2r}
\bar{\mathbf{x}}(t)=2\boldsymbol{\mathcal{J}}_n\mathbf{x}^r(t),
\end{equation*}
where  $\boldsymbol{\mathcal{J}}_n=\frac{1}{2} \boldsymbol{\mathcal{A}}\otimes \boldsymbol{I}_n$, and
$$\boldsymbol{\mathcal{A}}= \left[
\begin{array}{rrrr}
1 &  \imath &  \jmath   &   \kappa \\
1 &  -\imath &   \jmath  & -\kappa  \\
1 &  \imath &  -\jmath   &  -\kappa  \\
1 &  -\imath &  -\jmath   &   \kappa \\
\end{array}\right],$$
with $\boldsymbol{\mathcal{J}}_n^{\mathtt{H}}\boldsymbol{\mathcal{J}}_n=\boldsymbol{I}_{4n}$.

The { pseudo}-autocorrelation matrix of $\bar{\mathbf{x}}(t)$,  $\boldsymbol{ \Gamma}_{\bar{\mathbf{x}}}(t,s)= \E [\bar{\mathbf{x}}(t)\bar{\mathbf{x}}^\mathtt{H}(s)]$,  takes the form
\begin{equation*}
  \boldsymbol{ \Gamma}_{\bar{\mathbf{x}}}(t,s) =
  \begin{pmatrix}
  \boldsymbol{\Gamma}_{\mathbf{x}}(t,s) &   \boldsymbol{\Gamma}_{\mathbf{x}\mathbf{x}^{*}}(t,s) &     \boldsymbol{\Gamma}_{\mathbf{x}\mathbf{x}^{\imath}}(t,s)  &   \boldsymbol{\Gamma}_{\mathbf{x}\mathbf{x}^{\kappa}}(t,s) \\
  \boldsymbol{\Gamma}_{\mathbf{x}\mathbf{x}^{*}}^{*}(t,s) &    \boldsymbol{\Gamma}_{\mathbf{x}}^{*}(t,s)&   \boldsymbol{\Gamma}_{\mathbf{x}\mathbf{x}^{\kappa}}^{*}(t,s)  &   \boldsymbol{\Gamma}_{\mathbf{x}\mathbf{x}^{\imath}}^{*}(t,s) \\
  \boldsymbol{\Gamma}_{\mathbf{x}\mathbf{x}^{\imath}}^{\imath}(t,s) &    \boldsymbol{\Gamma}_{\mathbf{x}\mathbf{x}^{\kappa}}^{\imath}(t,s) &   \boldsymbol{\Gamma}_{\mathbf{x}}^{\imath}(t,s)  &   \boldsymbol{\Gamma}_{\mathbf{x}\mathbf{x}^{*}}^{\imath}(t,s) \\
  \boldsymbol{\Gamma}_{\mathbf{x}\mathbf{x}^{\kappa}}^{\kappa}(t,s) &    \boldsymbol{\Gamma}_{\mathbf{x}\mathbf{x}^{\imath}}^{\kappa}(t,s) &    \boldsymbol{\Gamma}_{\mathbf{x}\mathbf{x}^{*}}^{\kappa}(t,s)    &   \boldsymbol{\Gamma}_{\mathbf{x}}^{\kappa}(t,s) 
\end{pmatrix}
\end{equation*}
with  $\boldsymbol{ \Gamma}_{\mathbf{x}}(t,s)=\E [\mathbf{x}(t)\mathbf{x}^\mathtt{H}(s)]$, $\boldsymbol{ \Gamma}_{\mathbf{x} \mathbf{y}}(t,s)=\E [\mathbf{x}(t)\mathbf{y}^\mathtt{H}(s)]$, for any $\mathbf{x}(t),  \mathbf{y}(t) \in \mathbb{T}^{n}$.

Two key properties of tessarine random signals,  termed $\mathbb{T}_k$-properness ($k=1,2$),  have been established according to the vanishing of the  { pseudo}-correlation matrices $\boldsymbol{\Gamma}_{\mathbf{x}\mathbf{x}^{\nu}}(t,s)$, $\nu=*,\imath,\kappa$ \cite{Navarro2021}.
These properties are formally defined as follows.

\begin{defi}[$\mathbb{T}_k$-properness]\label{defTk}
$\mathbf{x}(t)\in \mathbb{T}^n$ is defined as $\mathbb{T}_1$-proper (or $\mathbb{T}_2$-proper) if, and only if,   the pseudo-correlations $\boldsymbol{\Gamma}_{\mathbf{x}\mathbf{x}^{\nu}}(t,s)=\boldsymbol{0}_{n\times n}$   for   $\nu=*,\imath,\kappa$ (respectively,  $\nu=\imath,\kappa$),  for all $t, s \in \mathbb{Z}$. 
Furthermore,   $\mathbf{x}(t)\in \mathbb{T}^{n_1}$ and $\mathbf{y}(t)\in \mathbb{T}^{n_2}$ are defined as cross  $\mathbb{T}_1$-proper (or cross $\mathbb{T}_2$-proper) if, and only if,  $\mathbf{\Gamma}_{\mathbf{x}\mathbf{y}^{\nu}}(t,s)=\boldsymbol{0}_{n\times n}$, for $\nu=*,\imath,\kappa$ (respectively, $\nu=\imath,\kappa$),  for all $t, s \in \mathbb{Z}$. 
In addition,  if $\mathbf{x}(t)$ and $\mathbf{y}(t)$ are both $\mathbb{T}_1$-proper  and cross $\mathbb{T}_1$-proper ( or $\mathbb{T}_2$-proper and  cross $\mathbb{T}_2$-proper),  then they are defined as  jointly $\mathbb{T}_1$-proper (respectively,  jointly $\mathbb{T}_2$-proper).
\end{defi}

\begin{remark}
The second-order statistical properties of a  $\mathbb{T}_k$-proper signal  $\mathbf{x}(t)\in \mathbb{T}^n$  can be fully characterized by either the signal itself (if $k=1$) or by the 2nD augmented vector 
$\mathbf{x}_2(t)=[\mathbf{x}^{\mathtt{T}}(t), \mathbf{x}^{\mathtt{H}}(t)]$ 
(if $k=2$). This results in a   reduction of the problem dimensionality  by a quarter  
or a half, respectively.  

This reduction in dimensionality directly impacts the complexity of optimal linear processing. 
While the optimal linear processing for general tessarine signals is WL processing, which operates on the full  4nD augmented vector $\bar{\mathbf{x}}(t)$,  under $\mathbb{T}_k$-proper properties ($k=1,2$),  the optimal linear processing simplifies to the $\mathbb{T}_k$-proper processing,  which means operating on a lower-dimensional vector  $\mathbf{x}_k(t)$,   where $\mathbf{x}_1(t)=\mathbf{x}(t)$ for $k=1$,  and $\mathbf{x}_2(t)=[\mathbf{x}^{\mathtt{T}}(t), \mathbf{x}^{\mathtt{H}}(t)]^{\mathtt{T}}$  for $k=2$.  This dimensionality reduction can result in significant computational savings. 

Consequently,   assessing
the $\mathbb{T}_k$-properness of a signal is key,  as it directly impacts the complexity of optimal linear processing.  
To empirically verify whether a signal satisfies $\mathbb{T}_1$- or $\mathbb{T}_2$-  proper conditions,  statistical tests have been developed and proposed in \cite{Navarro2021}.
\end{remark}

Turning our attention to signals characterized by  factorizable { pseudo}-autocorrelation functions,   this paper examines  the class of signals whose augmented { pseudo}-autocorrelation function is a factorizable kernel.  This family of signals,  referred to as {\it widely factorizable} signals,   is first introduced in a general context,  and then  particularized under $\mathbb{T}_k$-proper conditions.  

\begin{defi}[Widely factorizable signals]\label{def2}
$\mathbf{x}(t)\in \mathbb{T}^n$ is   widely factorizable if  its augmented {pseudo}-autocorrelation function $\boldsymbol{\Gamma}_{\bar{\mathbf{x}}}(t,s)$ can be expressed in the following factorized form:
\begin{equation}\label{fact4D}
\begin{array}{rl}
\boldsymbol{\Gamma}_{\bar{\mathbf{x}}}(t,s)= & \left\{
\begin{array}{lr}
\bar{\mathbf{A}}(t)\bar{\mathbf{B}}^{\mathtt{H}}(s),&  t\geq s\\
\bar{\mathbf{B}}(t)\bar{\mathbf{A}}^{\mathtt{H}}(s),&
t< s
 \end{array}
 \right.
\end{array}
\end{equation}
 \end{defi} 
with $\bar{\mathbf{A}}(t),  \bar{\mathbf{B}}(t) \in \mathbb{T}^{kn\times p}$.
  
\begin{remark}
Widely factorizable signals  represent a broad class that encompasses  both stationary and nonstationary signals, making them a versatile framework for both theoretical analysis and practical applications  \cite{Fernandez2019}.  
 Some important classes of signals that naturally exhibit widely factorizable correlation structures of the form \eqref{fact4D} are described below:
 \begin{enumerate}
 \item {\bf Linear dynamic models}. Consider $\mathbf{x}(t)\in \mathbb{T}^n $  described by a linear state model \cite{Jimenez2021}:
\begin{multline*}\label{statemodel}
\mathbf{x}(t+1)=\mathbf{F}_1(t)\mathbf{x}(t)+\mathbf{F}_2(t)\mathbf{x}^{*}(t)+\mathbf{F}_3(t)\mathbf{x}^{\imath}(t)
\\+\mathbf{F}_4(t)\mathbf{x}^{\kappa}(t)+\mathbf{w}(t),
\end{multline*}
where  $\mathbf{F}_1(t),  \mathbf{F}_2(t), \mathbf{F}_3(t), \mathbf{F}_4(t) \in \mathbb{T}^{n\times n}$ are known,   and $\mathbf{w}(t)$ is a tessarine noise.
Then, the augmented signal $\bar{\mathbf{x}}(t)$ can be expressed through the following WL state model:
$\bar{\mathbf{x}}(t+1)=\bar{\mathbf{F}}(t)\bar{\mathbf{x}}(t)+\bar{\mathbf{w}}(t)$, 
with
\begin{equation*}
\bar{\mathbf{F}}(t)=\left(
  \begin{array}{cccc}
    \mathbf{F}_1(t) & \mathbf{F}_2(t) & \mathbf{F}_3(t) & \mathbf{F}_4(t) \\
    \mathbf{F}_2^*(t) & \mathbf{F}_1^*(t) & \mathbf{F}_4^*(t) & \mathbf{F}_3^*(t) \\
    \mathbf{F}_3^{\imath}(t) & \mathbf{F}_4^{\imath}(t) & \mathbf{F}_1^{\imath}(t) & \mathbf{F}_2^{\imath}(t) \\
    \mathbf{F}_4^{\kappa}(t)& \mathbf{F}_3^{\kappa}(t) & \mathbf{F}_2^{\kappa}(t) & \mathbf{F}_1^{\kappa}(t) \\
  \end{array}
\right).
\end{equation*}
Thus,  $\mathbf{x}(t)$ is widely factorizable with $\bar{\mathbf{A}}(t)=\boldsymbol{\Phi}(t,0)$ and  $\bar{\mathbf{B}}(t)=\boldsymbol{\Gamma}_{\bar{\mathbf{x}}}(t)\boldsymbol{\Phi}^{-1}(t,0)$, where $\boldsymbol{\Phi}(t,0)=\prod\limits_{k=0}^{t-1}\bar{\mathbf{F}}(k)$, and $\boldsymbol{\Gamma}_{\bar{\mathbf{x}}}(t)=\boldsymbol{\Gamma}_{\bar{\mathbf{x}}}(t,t)$.

\item {\bf Wide-sense Markov processes}. Consider  ${\mathbf{x}}(t)\in \mathbb{T}^n $ whose augmented vector $\bar{\mathbf{x}}(t)$ is  a  wide-sense  Markov process of order $p\geq 1$ \cite{Navarro2021}.  Defining the  forward vector $\bar{\mathbf{z}}(t)=\left[\bar{\mathbf{x}}^{\mathtt{T}}(t), \bar{\mathbf{x}}^{\mathtt{T}}(t-1),\dots, \bar{\mathbf{x}}^{\mathtt{T}}(t-p+1)  \right]^{\mathtt{T}}$, the pseudo-autocorrelation function of $\bar{\mathbf{x}}(t)$    admits the factorization  \eqref{fact4D}, with $\bar{\mathbf{A}}(t)=\mathcal{V}\mathbf{G}(t)$ and  $\bar{\mathbf{B}}(t)=\mathcal{V}\boldsymbol{\Gamma}_{\bar{\mathbf{z}}}(t)\mathbf{G}^{-1}(t)$,  
where $\mathcal{V}=\left[ \boldsymbol{I}_{4n},\boldsymbol{0}_{4n\times 4n(p-1)} \right]$,  and $\mathbf{G}(t) = \prod\limits_{k=0}^{t-1}\boldsymbol{\Gamma}_{\bar{\mathbf{z}}}(k+1,k)\boldsymbol{\Gamma}^{-1}_{\bar{\mathbf{z}}}(k)$. 

\item {\bf Autoregressive moving average (ARMA) models.} Consider  ${\mathbf{x}}(t)\in \mathbb{T}^n $ an ARMA(p,q) process whose augmented vector $\bar{\mathbf{x}}(t)$ can be described through the following system:
$$\bar{\mathbf{x}}(t)=\sum\limits_{k=1}^p \bar{\mathbf{F}}_k\bar{\mathbf{x}}(t-k) + \sum\limits_{k=0}^q \bar{\mathbf{G}}_k\bar{\mathbf{w}}(t-k),$$
where   $ \bar{\mathbf{F}}_k,  \bar{\mathbf{G}}_k \in \mathbb{T}^{4n\times 4n}$ are deterministics,  $\bar{\mathbf{G}}_0=\boldsymbol{I}_{4}$,  and $\bar{\mathbf{w}}(t)$ is an augmented tessarine white noise vector.   Thus, $\mathbf{x}(t)$ is widely factorizable with $\bar{\mathbf{A}}(t)=\mathcal{V}\mathbf{H}^t$  and $\bar{\mathbf{B}}(t)=\mathcal{V}\boldsymbol{\Gamma}_{\bar{\mathbf{z}}}(t)\mathbf{H}^{-t}$,  
with   $\mathcal{V}=\left[ \boldsymbol{I}_{4n},\boldsymbol{0}_{4n\times 4n(p+q-1)} \right]$,  $\boldsymbol{\Gamma}_{\bar{\mathbf{z}}}(t)$ 
 the pseudo-autocorrelation function of    $\bar{\mathbf{z}}(t)=\left[\bar{\mathbf{x}}^{\mathtt{T}}(t), \bar{\mathbf{x}}^{\mathtt{T}}(t-1),\dots, \bar{\mathbf{x}}^{\mathtt{T}}(t-p+1), \bar{\mathbf{w}}^{\mathtt{T}}(t),\bar{\mathbf{w}}^{\mathtt{T}}(t-1),\right. $ $\left. \dots, \bar{\mathbf{w}}^{\mathtt{T}}(t-q+1)\right]^{\mathtt{T}}$,    
  $$\mathbf{H}=\left[\begin{array}{cccc} \bar{\mathbf{H}}_1 &  \boldsymbol{0}_{4np\times 4nq} \\ 
\boldsymbol{0}_{4nq\times 4np} &\bar{\mathbf{H}}_2
 \end{array}
\right],$$   with 
$$ \mathbf{H}_1=\left[\begin{array}{ccccc} \bar{\mathbf{F}}_1 & \bar{\mathbf{F}}_2 & \cdots & \bar{\mathbf{F}}_{p-1}& \bar{\mathbf{F}}_p \\
\boldsymbol{I}_{4n} &\boldsymbol{0}_{4n\times 4n} & \cdots & \boldsymbol{0}_{4n\times 4n}& \boldsymbol{0}_{4n\times 4n} \\
\boldsymbol{0}_{4n\times 4n} &\boldsymbol{I}_{4n} &\cdots & \boldsymbol{0}_{4n\times 4n} & \boldsymbol{0}_{4n\times 4n}\\
\vdots & \vdots & \ddots & \vdots  & \vdots \\
\boldsymbol{0}_{4n\times 4n} &\boldsymbol{0}_{4n\times 4n} &\cdots &  \boldsymbol{I}_{4n\times 4n} & \boldsymbol{0}_{4n\times 4n}
\end{array}
\right],
$$

$$ \mathbf{H}_2=\left[\begin{array}{ccccc} \bar{\mathbf{G}}_1 & \bar{\mathbf{G}}_2 & \cdots & \bar{\mathbf{G}}_{q-1}& \bar{\mathbf{G}}_q \\
 \boldsymbol{I}_{4n} &\boldsymbol{0}_{4n\times 4n} & \cdots & \boldsymbol{0}_{4n\times 4n}  & \boldsymbol{0}_{4n\times 4n}\\
 \boldsymbol{0}_{4n\times 4n} &\boldsymbol{I}_{4n} &\cdots & \boldsymbol{0}_{4n\times 4n} & \boldsymbol{0}_{4n\times 4n}\\
 \vdots & \vdots & \ddots & \vdots& \vdots \\
 \boldsymbol{0}_{4n\times 4n} &\boldsymbol{0}_{4n\times 4n} &\cdots & \boldsymbol{I}_{4n} & \boldsymbol{0}_{4n\times 4n}
\end{array}
\right].
$$

 \end{enumerate}
 \end{remark}
 

 The following result is easily deduced from Definition \ref{defTk} and Definition \ref{def2}.
 \begin{proposition}\label{prop1a}
If $\mathbf{x}(t)\in\mathbb{T}^n$  is $\mathbb{T}_k$-proper,  for $k=1,2$,  and widely factorizable then, the augmented {pseudo}-autocorrelation function of $\bar{\mathbf{x}}(t)$ takes the following form $\boldsymbol{\Gamma}_{\bar{\mathbf{x}}}(t,s)=\diag\left(\boldsymbol{\Gamma}_{\mathbf{x}_1}(t,s), \boldsymbol{\Gamma}_{\mathbf{x}_1}^{*}(t,s),  \boldsymbol{\Gamma}_{\mathbf{x}_1}^{\imath}(t,s), \boldsymbol{\Gamma}_{\mathbf{x}_1}^{\kappa}(t,s) \right)$, if $k=1$, or $\boldsymbol{\Gamma}_{\bar{\mathbf{x}}}(t,s)=\diag\left(\boldsymbol{\Gamma}_{\mathbf{x}_2}(t,s), \boldsymbol{\Gamma}_{\mathbf{x}_2}^{\mathtt{H}}(t,s) \right)$, if $k=2$,
with 
\begin{equation}\label{factk1}
\begin{array}{rl}
\boldsymbol{\Gamma}_{\mathbf{x}_k}(t,s)= & \left\{
\begin{array}{lr}
\mathbf{A}_k(t)\mathbf{B}_k^{\mathtt{H}}(s),&  t\geq s\\
\mathbf{B}_k(t)\mathbf{A}_k^{\mathtt{H}}(s),&
t< s
 \end{array}
 \right.
\end{array}
\end{equation}
where $\mathbf{A}_k(t)=\left[ \boldsymbol{I}_{kn}, \boldsymbol{0}_{kn\times (4-k)n}  \right]\bar{\mathbf{A}}(t)$ and $\mathbf{B}_k(t)=\left[ \boldsymbol{I}_{kn}, \boldsymbol{0}_{kn\times (4-k)n}  \right]\bar{\mathbf{B}}(t)$.
\end{proposition}

Proposition \ref{prop1a} motivates the introduction of  the class of $T_k$-proper factorizable signals.

\begin{defi}[$\mathbb{T}_k$-proper factorizable signals]\label{def3}  $\mathbf{x}(t)\in \mathbb{T}^n$ is   $\mathbb{T}_k$-proper  factorizable  if the { pseudo}-autocorrelation function of ${\mathbf{x}_k}(t)$,  $\boldsymbol{\Gamma}_{{\mathbf{x}_k}}(t,s)$,  can be expressed in the factorized form \eqref{factk1}.
 \end{defi}

          
          This section concludes with the definition of a specific product operation for tessarine vectors.

\begin{defi}
Given  $\mathbf{x}(t),\mathbf{y}(s)\in\mathbb{T}^n$,   $t,s\in \mathbb{Z}$,  the product $\star$ is defined as
\begin{multline*}\label{star product}
 \mathbf{x}(t) \star  \mathbf{y}(s)=\mathbf{x}_r(t) \circ \mathbf{y}_r(s)+\imath\mathbf{x}_{\imath}(t)\circ \mathbf{y}_{\imath}(s)+\jmath \mathbf{x}_{\jmath}(t)\circ \mathbf{y}_{\jmath}(s)\\ +\kappa \mathbf{x}_{\kappa}(t)\circ \mathbf{y}_{\kappa}(s).
 \end{multline*} 
\end{defi}

\begin{property}\label{propstar}
Given  $\mathbf{x}(t),\mathbf{y}(s)\in\mathbb{T}^n$,  $t,s\in \mathbb{Z}$,  the augmented vector of $\mathbf{x}(t)\star\mathbf{y}(s)$ is given by $\overline{\mathbf{x}(t)\star\mathbf{y}(s)}=\boldsymbol{\mathcal{D}}^{{\mathbf{x}}}(t)\bar{\mathbf{y}}(s)$,   where $\boldsymbol{\mathcal{D}}^{{\mathbf{x}}}(t)=\boldsymbol{\mathcal{J}}_n\diag\left(\mathbf{x}^r(t)\right)\boldsymbol{\mathcal{J}}_n^{\texttt{H}}$.
\end{property}

\section{Problem statement}\label{formulation}

This section formulates the problem of estimating a $\mathbb{T}_k$-proper factorizable signal,  based on fading measurements received from $R$ sensors.

Consider a $\mathbb{T}_k$-proper  factorizable  tessarine random signal $\{ \mathbf{x}(t)\in \mathbb{T}^n,  \  t\in \mathbb{Z}\}$, observed by $R$ sensors.   Each sensor detects the signal with added noise and attenuation due to fading,  as described in the  model below:
\begin{equation}\label{obs}
\mathbf{y}^{(\alpha)}(t)=   \boldsymbol{\gamma}^{(\alpha)}(t) \star\mathbf{x}(t)+\mathbf{v}^{(\alpha)}(t), \ \ t\geq 1,  
\end{equation}
where,  at each sensor $ \alpha=1,\dots,R$,   $ \boldsymbol{\gamma}^{(\alpha)}(t)=[\gamma_{1}^{(\alpha)}(t), \dots,\gamma_{n}^{(\alpha)}(t) ]^{\texttt{T}}\in \mathbb{T}^{n}$  represents the fading gain vector of the $\alpha$-th sensor, and $\mathbf{v}^{(\alpha)}(t)\in \mathbb{T}^{n}$  is a tessarine white noise with {pseudo}-variance $\mathbf{R}^{(\alpha)}(t)$.
The following hypotheses  are assumed on these processes:
\begin{enumerate}

\item[H.1.] $\gamma_{j}^{(\alpha)}(t)=\gamma_{j,\mathrm{r}}^{(\alpha)}(t)+\imath \gamma_{j,\imath}^{(\alpha)}(t)+\jmath \gamma_{j,\jmath}^{(\alpha)}(t)+\kappa \gamma_{j,\kappa}^{(\alpha)}(t)$, for $j=1,\ldots, n$, where  $\gamma_{j,\nu}^{(\alpha)}(t)$, $\nu=\mathrm{r}, \imath,\jmath,\kappa$, are  independent real random variables that  take values over any subinterval of $[0,1]$  according to an arbitrary distribution with known second-order statistics.

\item[H.2.] $\boldsymbol{\gamma}^{(\alpha)}(t)$, $\boldsymbol{\gamma}^{(\beta)}(s)$ are independent for any  $\alpha\neq \beta$ or  $t\neq s$.


\item[H.3.] $\boldsymbol{\gamma}^{(\alpha)}(t)$, $\mathbf{x}(t)$, and $\mathbf{v}^{(\alpha)}(t)$ are mutually  independent.

\item[H.4.] 
$\boldsymbol{\Gamma}_{\mathbf{v}}^{(\alpha \beta)}(t,s)=\E[\mathbf{v}^{(\alpha)}(t)\mathbf{v}^{(\beta)^\texttt{H}}(s)]=\mathbf{R}^{(\alpha \beta)}(t)\delta_{t,s}$,  for any $\alpha,\beta=1,\dots,R$.

\item[H.5.] $\mathbf{x}(t)$ and $\mathbf{y}^{(\alpha)}(t)$  are  jointly $\mathbb{T}_k$-proper,  for  each sensor $\alpha=1,\dots,R$.
\end{enumerate}

The formulated problem is closely related to practical scenarios in multisensor networks operating under unreliable conditions.  Specifically,  the fading gain vectors, $\boldsymbol{\gamma}^{(\alpha)}(t)$,   capture the random attenuation and fluctuations commonly encountered in wireless and underwater communication channels, where measurements may be distorted by multipath propagation, scattering, or medium inhomogeneities.

\begin{remark}\label{Rem3}
Let $\E[\gamma_{j,\nu}^{(\alpha)}(t)]=\mu_{j,\nu}^{(\alpha)}(t)$ and $\var[\gamma_{j,\nu}^{(\alpha)}(t)]=\sigma_{j,\nu}^{(\alpha)^2}(t)$,  for $j=1,\dots,n$,  $\nu=\mathrm{r}, \imath,\jmath,\kappa$.  From hypotheses H.1-H.2, it can be readily  verified that $\E[\gamma_{j_1,\nu_1}^{(\alpha_1)}(t)\gamma_{j_2,\nu_2}^{(\alpha_2)}(s)]=\sigma_{j_1,\nu_1}^{(\alpha_1)^2}(t)+\mu_{j_1,\nu_1}^{(\alpha_1)^2}(t)$,  for $j_1=j_2, \alpha_1=\alpha_2, \nu_1=\nu_2, t=s$, and $\E[\gamma_{j_1,\nu_1}^{(\alpha_1)}(t)\gamma_{j_2,\nu_2}^{(\alpha_2)}(s)]= \mu_{j_1,\nu_1}^{(\alpha_1)}(t)\mu_{j_2,\nu_2}^{(\alpha_2)}(s)$, otherwise.
\end{remark}

The Proposition below provides  a characterization of the joint $\mathbb{T}_k$-properness of $\mathbf{x}(t)$ and $\mathbf{y}^{(\alpha)}(t)$ in terms of second-order statistics of $\gamma_{j,\nu}^{(\alpha)}(t)$.

\begin{proposition}\label{prop1}
Let $\mathbf{x}(t)\in \mathbb{T}^n$ be an $n$D signal  observed  at each sensor $\alpha$ according to the equation \eqref{obs}. Then,
\begin{enumerate}
\item  $\mathbf{x}(t)$, $\mathbf{y}^{(\alpha)}(t)$ are jointly $\mathbb{T}_1$-proper if, and only if, $\mathbf{x}(t)$ and $\mathbf{v}^{(\alpha)}(t)$ are $\mathbb{T}_1$-proper, and $\mu_{j,\nu}^{(\alpha)}(t)=\mu_{j}^{(\alpha)}(t)$, $\sigma_{j,\nu}^{(\alpha)}(t)=\sigma_{j}^{(\alpha)}(t)$, $\forall t\in \mathbb{Z}, j=1,\dots,n$,  $\alpha=1,\dots,R$, and $\nu=\mathrm{r}, \imath, \jmath, \kappa$. 
\item  $\mathbf{x}(t)$, $\mathbf{y}^{(\alpha)}(t)$ are jointly $\mathbb{T}_2$-proper if, and only if, $\mathbf{x}(t)$ and $\mathbf{v}^{(\alpha)}(t)$ are $\mathbb{T}_2$-proper, and $\mu_{j,\mathrm{r}}^{(\alpha)}(t)=\mu_{j,\jmath}^{(\alpha)}(t)$, $\mu_{j,\imath}^{(\alpha)}(t)=\mu_{j,\kappa}^{(\alpha)}(t)$, $\sigma_{j,\mathrm{r}}^{(\alpha)}(t)=\sigma_{j,\jmath}^{(\alpha)}(t)$, $\sigma_{j,\imath}^{(\alpha)}(t)=\sigma_{j,\kappa}^{(\alpha)}(t)$, $\forall t\in \mathbb{Z}, j=1,\dots,n$, and $\alpha=1,\dots,R$.  
\end{enumerate}
\end{proposition}

In this setting, the goal is to estimate the signal $\mathbf{x}(t)$ from the observations $\{ \mathbf{y}^{(\alpha)}(1), \dots,  \mathbf{y}^{(\alpha)}(s)  \}$, $\alpha=1,\dots, R$.  To achieve this,  a $\mathbb{T}_k$-proper processing technique is  employed to devise 
a distributed fusion estimation algorithm based on the second-order statistical information of the processes involved.

In general,   the $\mathbb{T}_k$-proper processing methodology  requires operating on the observations \cite{Jimenez2021, Fernandez2023}:
\begin{equation}\label{yk}
\mathbf{y}_k^{(\alpha)}(t)=
\mathcal{C}_k \boldsymbol{\mathcal{D}}^{\boldsymbol{\gamma}^{(\alpha)}}(t)\bar{\mathbf{x}}(t)
+\mathbf{v}_k^{(\alpha)}(t),  \quad t\geq 1,  
\end{equation}
for $\alpha=1,\dots, R$,  $k=1,2$,  with $\mathcal{C}_k=\left[ \boldsymbol{I}_{kn},\boldsymbol{0}_{kn\times (4-k)n} \right]$,  and $\boldsymbol{\mathcal{D}}^{\boldsymbol{\gamma}^{(\alpha)}}(t)$ defined in Property \ref{prop2}. 

Although $\mathbf{y}_k^{(\alpha)}(t)$ is a $kn$D vector ($k=1,2$),  
 equation  \eqref{yk} requires working with the $4n$D augmented signal vector  $\bar{\mathbf{x}}(t)$. As an alternative approach,  a 
process $\mathbf{z}_k^{(\alpha)}(t)$ is introduced, characterized by a linear equation that involves lower-dimensional processes while preserving the second-order statistical properties of 
$\mathbf{y}_k^{(\alpha)}(t)$.  Specifically,  consider the equation
\begin{equation}\label{zk}
\mathbf{z}_k^{(\alpha)}(t)=\mathbf{H}_k^{(\alpha)}(t)\mathbf{x}_k(t)+\mathbf{w}_k^{(\alpha)}(t), \quad t\geq 1,  
\end{equation}
for $\alpha=1,\dots, R$,  $k=1,2$,  where  $\mathbf{w}_k^{(\alpha)}(t)$ is a white noise independent of $\mathbf{x}_k(t)$, and 
$\boldsymbol{\Gamma}_{\mathbf{w}_k}^{(\alpha\beta)}(t)=
\mathbf{R}_k^{(\alpha\beta)}(t)+\boldsymbol{\Sigma}_k^{(\alpha)}(t) \delta_{\alpha,\beta}$. Moreover, 

\begin{itemize}
\item For $k=1$, 
\begin{itemize}
\item $\mathbf{H}_1^{(\alpha)}(t)=\diag\left(\boldsymbol{\mu}^{(\alpha)}(t)\right)$, where $\boldsymbol{\mu}^{(\alpha)}(t)=[\mu_1^{(\alpha)}(t),\ldots,\mu_n^{(\alpha)}(t)]^{\texttt{T}}$, with   $\mu_j^{(\alpha)}(t)$ as in Proposition \ref{prop1}.
\item $\mathbf{R}_1^{(\alpha\beta)}(t)=\boldsymbol{\Gamma}_{\mathbf{v}}^{(\alpha\beta)}(t)$, and 
 $\boldsymbol{\Sigma}_1^{(\alpha)}(t)=\diag(\boldsymbol{\theta}^{(\alpha)}(t))$, where $\boldsymbol{\theta}^{(\alpha)}(t)=[\theta_1^{(\alpha)}(t),\ldots,\theta_n^{(\alpha)}(t) ]^{\texttt{T}}$, with   $\theta_j^{(\alpha)}(t)=4\sigma_{j,\mathrm{r}}^{(\alpha)^2}(t)\E\left[x_{j,\mathrm{r}}^2(t)\right]$,   $j=1,\dots,n$.
\end{itemize}

\item For $k=2$,
\begin{itemize}
\item $\mathbf{H}_2^{(\alpha)}(t)=\left[h_{ij}^{(\alpha)}(t)\right]$, where 
\begin{equation*}
\begin{split}
h_{ii}^{(\alpha)}(t)= & h_{i+n,i+n}^{(\alpha)}(t)=\frac{1}{2}\left(\mu_{i,\mathrm{r}}^{(\alpha)}(t)+\mu_{i,\imath}^{(\alpha)}(t)\right),\\
h_{i,i+n}^{(\alpha)}(t)= & h_{i+n,i}^{(\alpha)}(t)=\frac{1}{2}\left(\mu_{i,\mathrm{r}}^{(\alpha)}(t)-\mu_{i,\imath}^{(\alpha)}(t)\right),\\
h_{i,j}^{(\alpha)}(t)= & 0,  \text{ for the rest.}
\end{split}
\end{equation*}
for $  i=1,\dots,n$,  $\alpha=1,\dots,R$,.
\item $\mathbf{R}_2^{(\alpha\beta)}(t)=\boldsymbol{\Gamma}_{\mathbf{v}_2}^{(\alpha\beta)}(t)$, with $\mathbf{v}_2^{(\alpha)}(t)=[\mathbf{v}^{(\alpha)^{\texttt{T}}}(t),\mathbf{v}^{(\alpha)^{\texttt{H}}}(t) ]^{\texttt{T}}$, and 
 $\boldsymbol{\Sigma}_2^{(\alpha)}(t)=\left[\lambda_{ij}^{(\alpha)}(t)\right]$, where 
\begin{equation*}
\begin{split}
\lambda_{ii}^{(\alpha)}(t)= & \lambda_{i+n,i+n}^{(\alpha)}(t)=\phi_{i,\mathrm{r}}^{(\alpha)}(t)+\phi_{i,\imath}^{(\alpha)}(t),\\
\lambda_{i,i+n}^{(\alpha)}(t)= & \lambda_{i+n,i}^{(\alpha)}(t)=\phi_{i,\mathrm{r}}^{(\alpha)}(t)-\phi_{i,\imath}^{(\alpha)}(t),\\
\lambda_{ij}^{(\alpha)}(t)= & 0, \text{ for the rest},
\end{split}
\end{equation*}
for $  i=1,\dots,n$,  $\alpha=1,\dots,R$,  with $\phi_{j,\nu}^{(\alpha)}(t)=2 \sigma_{j,\nu}^{(\alpha)^2}(t)E\left[x_{j,\nu}^2(t)\right]$, for  $j=1,\dots,n$,\ \ $\nu=\mathrm{r}, \imath$.
\end{itemize}
\end{itemize}

The following Proposition establishes the equivalence between the {pseudo}-autocorrelation functions of $\mathbf{y}_k^{(\alpha)}(t)$ and $\mathbf{z}_k^{(\alpha)}(t)$, as well as the pseudo-cross-correlation function of these processes and $x_k(t)$.
\begin{proposition}\label{prop2}
For any two sensors $\alpha, \beta=1\dots, R$,  the following relations are satisfied: 
\begin{equation*}
\begin{split}
\boldsymbol{\Gamma}_{\mathbf{y}_k\mathbf{x}_k}^{(\alpha)}(t,s)= & \boldsymbol{\Gamma}_{\mathbf{z}_k\mathbf{x}_k}^{(\alpha)}(t,s),\\
\boldsymbol{\Gamma}_{\mathbf{y}_k}^{(\alpha\beta)}(t,s)= & \boldsymbol{\Gamma}_{\mathbf{z}_k}^{(\alpha\beta)}(t,s).
\end{split}
\end{equation*}

\end{proposition}

\section{$\mathbb{T}_k$-proper distributed fusion methodology}\label{distributed section}

In the conditions established in Section \ref{formulation}, this section focuses on  the distributed fusion linear MMSE estimation problem based on a $\mathbb{T}_k$-proper processing.

There are two phases in the distributed fusion method: First,  $\mathbb{T}_k$-proper local linear MMSE estimation algorithms for $\mathbf{x}_k(t)$ are derived by using the observed data from each sensor.  Then, 
a global  distributed fusion estimator is obtained by combining the local estimators through a fusion criterion that utilizes weighted matrices in the MMSE sense.  These two phases are further detailed in the following subsections.

\subsection{$\mathbb{T}_k$-proper local estimation algorithms}\label{local section}

For each sensor $\alpha=1,\ldots,R$,  this section introduces recursive algorithms for computing the  linear MMSE estimators  $\hat{{\mathbf{x}}}_k^{(\alpha)}(t|s)$ of $\mathbf{x}_k(t)$,  based on the set of observations $\{\mathbf{y}_k^{(\alpha)}(1), \dots , \mathbf{y}_k^{(\alpha)}(s)\}$,  for $k = 1, 2$,  and their corresponding  error pseudo-covariance matrices $\mathbf{P}_k^{(\alpha)}(t|s)$.   The mathematical formulations for these algorithms, covering filtering ($t=s$), prediction ($t>s$), and smoothing ($t<s$) scenarios are summarized in the following three Theorems,  respectively.  

\begin{remark}
The $\mathbb{T}_k$-proper  local linear MMSE estimators  of $\mathbf{x}(t)$ 
are extracted from the first $n$ components of the global estimators $\hat{{\mathbf{x}}}_k^{(\alpha)}(t|s)$.  Similarly,  the associated  $\mathbb{T}_k$-proper error pseudo-covariances,  ${\mathbf{P}}^{(\alpha)^{\mathbb{T}_k}}(t|s)$,  are determined by the first $n\times n$  elements of ${\mathbf{P}}_k^{(\alpha)}(t|s)$.
\end{remark}
 
 \begin{theorem}[Local  linear MMSE filter]\label{ThLocalFilter}
For the $\alpha$-th sensor,  with $\alpha=1,\ldots,R$,  the   linear MMSE filter  $\hat{{\mathbf{x}}}_k^{(\alpha)}(t|t)$
and its error  pseudo-variance matrix ${\mathbf{P}}_k^{(\alpha)}(t|t)$ are computed as follows:
\begin{align}\label{local filter}  
  \hat{{\mathbf{x}}}_k^{(\alpha)}(t|t) & = \mathbf{A}_k(t)\mathbf{e}^{(\alpha)}(t),\quad t\geq 1,
  \\  \label{P(t|t) local filter}
  {\mathbf{P}}_k^{(\alpha)}(t|t) &=\mathbf{A}_k(t) \left[\mathbf{B}_k^{\texttt{H}}(t)-\mathbf{Q}^{(\alpha)}(t)\mathbf{A}_k^{\texttt{H}}(t) \right],\quad t\geq 1,
\end{align}
where  $ \mathbf{e}^{(\alpha)}(t)$ and $ \mathbf{Q}^{(\alpha)}(t)$ are computed  recursively as follows:
\begin{align}
   \label{ek}
   \mathbf{e}^{(\alpha)}(t) & =  \mathbf{e}^{(\alpha)}(t-1)+\mathbf{J}_k^{(\alpha)}(t)\boldsymbol{\Omega}_k^{(\alpha)^{-1}}(t)\boldsymbol{\varepsilon}_k^{(\alpha)}(t),
 \\ \label{Innok}
  \boldsymbol{\varepsilon}_k^{(\alpha)}(t) & =\mathbf{y}_k^{(\alpha)}(t)-\mathbf{H}_k^{(\alpha)}(t)\mathbf{A}_k(t)\mathbf{e}^{(\alpha)}(t-1),
 \\ \label{Jk}
  \mathbf{J}_k^{(\alpha)}(t) & =\left[\mathbf{B}_k^{\texttt{H}}(t)-\mathbf{Q}^{(\alpha)}(t-1)\mathbf{A}_k^{\texttt{H}}(t) \right]\mathbf{H}_k^{(\alpha)^{\texttt{H}}}(t),
\\ \label{Qk}
  \mathbf{Q}^{(\alpha)}(t) & =  \mathbf{Q}^{(\alpha)}(t-1)+\mathbf{J}_k^{(\alpha)}(t)\boldsymbol{\Omega}_k^{(\alpha)^{-1}}(t)\mathbf{J}_k^{(\alpha)^{\texttt{H}}}(t),
 \\ \label{Omega}
   \boldsymbol{\Omega}_k^{(\alpha)}(t) & = \mathbf{R}_k^{(\alpha)}(t)+ \boldsymbol{\Sigma}_k^{(\alpha)}(t)+\mathbf{H}_k^{(\alpha)}(t)\mathbf{A}_k(t)  \mathbf{J}_k^{(\alpha)}(t)
 \end{align}
for $t\geq 1$,  with initial conditions $\mathbf{e}^{(\alpha)}(0)=\boldsymbol{0}_{p}$  and $\mathbf{Q}^{(\alpha)}(0)=\boldsymbol{0}_{p\times p}$.  Moreover, the matrices $\mathbf{H}_k^{(\alpha)}(t)$, $\mathbf{R}_k^{(\alpha)}(t)$, and $\boldsymbol{\Sigma}_k^{(\alpha)}(t)$ are as defined in Equation \eqref{zk}.

\end{theorem}

\begin{theorem}[Local linear  MMSE predictor]\label{LocalPredictor}
For the $\alpha$-th sensor,  $\alpha=1,\ldots,R$,  the   linear MMSE predictor  $\hat{{\mathbf{x}}}_k^{(\alpha)}(t|s)$, with $t>s$, 
and its error  pseudo-variance matrix ${\mathbf{P}}_k^{(\alpha)}(t|s)$ are computed as follows:
\begin{align}\label{local predictor}
    \hat{{\mathbf{x}}}_k^{(\alpha)}(t|s) & = \mathbf{A}_k(t)\mathbf{e}^{(\alpha)}(s),
    \\  \label{P(t|N) local predictor}
    {\mathbf{P}}_k^{(\alpha)}(t|s)& =\mathbf{A}_k(t) \left[\mathbf{B}_k^{\texttt{H}}(t)-\mathbf{Q}^{(\alpha)}(s)\mathbf{A}_k^{\texttt{H}}(t) \right],
  \end{align}
with $\mathbf{e}^{(\alpha)}(s)$ and $\mathbf{Q}^{(\alpha)}(s)$  recursively computed   from \eqref{ek} and \eqref{Qk}, respectively.

\end{theorem}

 \begin{theorem}[Local  linear MMSE smoother]\label{LocalSmoother}
For the $\alpha$-th sensor,  $\alpha=1,\ldots,R$,   the linear MMSE smoother, $\hat{{\mathbf{x}}}_k^{(\alpha)}(t|s)$, with  $t<s$,  and its error pseudo-variance matrix ${\mathbf{P}}_k^{(\alpha)}(t|s)$ are computed as follows:
\begin{align}\label{local smoother}
   \hat{{\mathbf{x}}}_k^{(\alpha)}(t|s) & = \hat{{\mathbf{x}}}_k^{(\alpha)}(t|s-1) + 
\mathbf{L}_k^{(\alpha)}(t,s)\boldsymbol{\varepsilon}_k^{(\alpha)}(s),
\\ \label{P(t|N) local smoother}
  {\mathbf{P}}_k^{(\alpha)}(t|s) &= {\mathbf{P}}_k^{(\alpha)}(t|s-1)-\mathbf{L}_k^{(\alpha)}(t,s) \boldsymbol{\Omega}_k^{(\alpha)}(s) \mathbf{L}_k^{(\alpha)^{\texttt{H}}}(t,s),
\end{align}
with initial conditions the linear MMSE filter $ \hat{{\mathbf{x}}}_k^{(\alpha)}(t|t) $ and
its error pseudo-variance  ${\mathbf{P}}_k^{(\alpha)}(t|t)$ given by  \eqref{local filter} and  \eqref{P(t|t) local filter}, respectively.  Moreover,   $\boldsymbol{\varepsilon}_k^{(\alpha)}(s)$ is recursively computed from \eqref{Innok}, and the smoothing gain matrix $ \mathbf{L}_k^{(\alpha)}(t,s) $ is calculated as follows:
\begin{multline}\label{Lk}
  \mathbf{L}_k^{(\alpha)}(t,s)  =\left[\mathbf{B}_k(t)-\mathbf{M}_k^{(\alpha)}(t,s-1) \right] \\ \times
   \mathbf{A}_k^{\texttt{H}}(s)\mathbf{H}_k^{(\alpha)^{\texttt{H}}}(s)\boldsymbol{\Omega}_k^{(\alpha)^{-1}}(s),
    \end{multline}
    \begin{equation}\label{Mk}
  \mathbf{M}_k^{(\alpha)}(t,s)  =  \mathbf{M}_k^{(\alpha)}(t,s-1)+\mathbf{L}_k^{(\alpha)}(t,s)\mathbf{J}_k^{(\alpha)^{\texttt{H}}}(s), 
 \end{equation}
for $t<s$,  with initial condition $ \mathbf{M}_k^{(\alpha)}(t,t)=\mathbf{A}_k(t)\mathbf{Q}^{(\alpha)}(t)$, where $\mathbf{Q}^{(\alpha)}(t)$ is obtained from the equation \eqref{Qk}. Moreover,  the matrices   $\mathbf{J}_k^{(\alpha)}(s)$  and $ \boldsymbol{\Omega}_k^{(\alpha)}(s)$ are computed from  \eqref{Jk} and \eqref{Omega},   respectively.

\end{theorem}

\subsection{$\mathbb{T}_k$-proper  distributed estimation algorithms}\label{DistributedAlgorithms}
In a second stage, local linear MMSE estimators $\hat{{\mathbf{x}}}_k^{(\alpha)}(t|s)$, $\alpha=1,\dots, R$, are fused  to produce a  global distributed estimator $\hat{\mathbf{x}}_k^D{(t|s)}$. Then, the first $n$-components of $\hat{\mathbf{x}}_k^D{(t|s)}$ are extracted to generate the desired $\mathbb{T}_k$-proper distributed estimator $\hat{\mathbf{x}}^D{(t|s)}$. Specifically,  using a matrix weighted fusion criterion,  $\hat{\mathbf{x}}_k^D{(t|s)}$ is given by the expression

\begin{equation}\label{Distest}
 \hat{{\mathbf{x}}}_k^{D}(t|s)=\boldsymbol{\mathcal{H}}_k(t,s) \vec{\hat{\mathbf{x}}}_k(t|s),
\end{equation}
where $\vec{\hat{\mathbf{x}}}_k(t|s)=\left[{\hat{\mathbf{x}}}_k^{(1)^{\texttt{T}}}(t|s),\ldots,\hat{{\mathbf{x}}}_k^{(R)^{\texttt{T}}}(t|s)\right]^{\texttt{T}}$,  $\boldsymbol{\mathcal{H}}_k(t,s)= \boldsymbol{\mathcal{O}}_k(t,s)  \boldsymbol{\mathcal{V}}_k^{-1}(t,s)$, with 
\begin{equation*}\label{Ok}
  \begin{aligned}
\boldsymbol{\mathcal{O}}_k(t,s) &= E[{\mathbf{x}}_k(t)\vec{\hat{\mathbf{x}}}_k^{{\texttt{H}}}(t|s)]=\left[\boldsymbol{\mathcal{V}}_k^{(11)}(t,s),\ldots,\boldsymbol{\mathcal{V}}_k^{(RR)}(t,s)\right],\\
\boldsymbol{\mathcal{V}}_k(t,s)&=  E[\vec{\hat{\mathbf{x}}}_k(t|s)\vec{\hat{\mathbf{x}}}_k^{{\texttt{H}}}(t|s)] =\left[\boldsymbol{\mathcal{V}}_k^{(\alpha\beta)}(t,s)\right]_{\alpha,\beta=1,\ldots,R},
  \end{aligned}
\end{equation*}
and $\boldsymbol{\mathcal{V}}_k^{(\alpha\beta)}(t,s)=E[\hat{{\mathbf{x}}}_k^{(\alpha)}(t|s)\hat{{\mathbf{x}}}_k^{(\beta)^{\texttt{H}}}(t|s)]$.     

Furthermore, the corresponding error pseudo-variance matrix $\mathbf{P}^D(t|s)$ is given by the first $n\times n$-submatrix of $ {\mathbf{P}}_k^D(t|s)$, where
\begin{equation}\label{Disterror}
  {\mathbf{P}}_k^{D}(t|s)=
  \mathbf{A}_k(t)\mathbf{B}_k^{\texttt{H}}(t)-\boldsymbol{\mathcal{O}}_k(t,s) \boldsymbol{\mathcal{V}}_k^{-1}(t,s)\boldsymbol{\mathcal{O}}_k^{\texttt{H}}(t,s).
\end{equation}

Expressions \eqref{Distest} and \eqref{Disterror} depend on the pseudo-cross-covariance matrices 
$\boldsymbol{\mathcal{V}}_k^{(\alpha\beta)}(t,s)$.  Recursive algorithms for their computation are detailed in Theorems \ref{teo distributed filter} - \ref{teo distributed smoother}.

\begin{theorem}[Filtering pseudo-cross-covariance matrices]\label{teo distributed filter}
The pseudo-cross-covariance $\boldsymbol{\mathcal{V}}_k^{(\alpha\beta)}(t)=\boldsymbol{\mathcal{V}}_k^{(\alpha\beta)}(t)$ 
is given by
\begin{equation}\label{Kij(t,t)}
 \boldsymbol{\mathcal{V}}_k^{(\alpha\beta)}(t)=\mathbf{A}_k(t)\mathbf{Q}^{(\alpha\beta)}(t)\mathbf{A}_k^{\texttt{H}}(t),\quad t\geq 1,
\end{equation}
where  
\begin{multline}\label{Qij(t)}
 \mathbf{Q}^{(\alpha\beta)}(t)=  \mathbf{Q}^{(\alpha\beta)}(t-1)+\mathbf{J}_k^{(\alpha)}(t)\boldsymbol{\Omega}_k^{(\alpha)^{-1}}(t)\mathbf{J}_k^{(\beta\alpha)^{\texttt{H}}}(t)\\+\mathbf{J}_k^{(\alpha\beta)}(t-1,t)\boldsymbol{\Omega}_k^{(\beta)^{-1}}(t)\mathbf{J}_k^{(\beta)^{\texttt{H}}}(t),\quad t\geq 1,
\end{multline}
with initial condition $\mathbf{Q}^{(\alpha\beta)}(0)=\boldsymbol{0}_{p \times p}$, and where
\begin{equation}\label{J^ij(t-1,t)}
 \mathbf{J}_k^{(\alpha\beta)}(t-1,t) =  \left[  \mathbf{Q}^{(\alpha)}(t-1) - \mathbf{Q}^{(\alpha\beta)}(t-1)\right]\mathbf{A}_k^{\texttt{H}}(t)\mathbf{H}_k^{(\beta)^{\texttt{H}}}(t), 
 \end{equation}
 for $t\geq 2$,  with initial condition $ \mathbf{J}_k^{(\alpha\beta)}(0,1)=\boldsymbol{0}_{p\times kn}$,  and 
 \begin{equation}\label{J^ij(t)}
 \mathbf{J}_k^{(\alpha\beta)}(t) =  \mathbf{J}_k^{(\alpha\beta)}(t-1,t)+ \mathbf{J}_k^{(\alpha)}(t)\boldsymbol{\Omega}_k^{(\alpha)^{-1}}(t)\boldsymbol{\Omega}_k^{(\alpha\beta)}(t),
\end{equation}
 for $t\geq 1$,   with
\begin{multline}\label{Omega^ij(t)}
   \boldsymbol{\Omega}_k^{(\alpha\beta)}(t) =  \mathbf{R}_k^{(\alpha \beta)}(t)+ \boldsymbol{\Sigma}_k^{(\alpha)}(t) \delta_{\alpha,\beta}\\ +\mathbf{H}_k^{(\alpha)}(t)\mathbf{A}_k(t)\left[\mathbf{J}_k^{(\beta)}(t) - \mathbf{J}_k^{(\alpha\beta)}(t-1,t) \right], 
\end{multline}
for $ t\geq 2$,  and  $\boldsymbol{\Omega}_k^{(\alpha\beta)}(1)=  \mathbf{R}_k^{(\alpha \beta)}(1)+ \boldsymbol{\Sigma}_k^{(\alpha)}(1) \delta_{\alpha,\beta}+\mathbf{H}_k^{(\alpha)}(1)\mathbf{A}_k(1)\mathbf{B}_k^{\texttt{H}}(1) \mathbf{H}_k^{(\beta)^\texttt{H}}(1)$.
\end{theorem}

\begin{remark}
For $\alpha=\beta$,   $\mathbf{J}_k^{(\alpha\beta)}(t)= \mathbf{J}_k^{(\alpha)}(t) $, $\mathbf{Q}^{(\alpha\beta)}(t)=\mathbf{Q}^{(\alpha)}(t)$ and $\boldsymbol{\Omega}_k^{(\alpha\beta)}(t)=\boldsymbol{\Omega}_k^{(\alpha)}(t)$ , which are computed from \eqref{Jk},  \eqref{Qk} and \eqref{Omega}, respectively.  Moreover,  it is trivial  that,  for $\alpha=\beta$,  $ \mathbf{J}_k^{(\alpha\beta)}(t-1,t)=\mathbf{J}_k^{(\alpha)}(t-1,t)=\boldsymbol{0}_{p\times kn}$.
\end{remark}

\begin{theorem}[Prediction pseudo-cross-covariance matrices]\label{teo distributed predictor}
The  pseudo-cross-covariances  $\boldsymbol{\mathcal{V}}_k^{(\alpha\beta)}(t,s)$, for $t>s$,  are given by
\begin{equation*}\label{K_k^{(ij)}(t,s)}
\boldsymbol{\mathcal{V}}_k^{(\alpha\beta)}(t,s)= \mathbf{A}_k(t)\mathbf{Q}^{(\alpha\beta)}(s) \mathbf{A}_k^\texttt{H}(t),  \quad t>s,
\end{equation*}
where $\mathbf{Q}^{(\alpha\beta)}(s)$ is computed as in Theorem \ref{teo distributed filter}.    
\end{theorem}

\begin{theorem}[Smoothing  pseudo-cross-covariance matrices]\label{teo distributed smoother}
The pseudo-cross-covariances $\boldsymbol{\mathcal{V}}_k^{(\alpha\beta)}(t,s)$, for $t<s$, are given by
\begin{multline}\label{K_k^(ij) smooth}
\boldsymbol{\mathcal{V}}_k^{(\alpha \beta)}(t,s)= \boldsymbol{\mathcal{V}}_k^{(\alpha \beta)}(t,s-1)+\mathbf{L}_k^{(\alpha)}(t,s)\boldsymbol{\mathcal{L}}_k^{(\beta \alpha)^{\texttt{H}}}(t,s)\\ +\boldsymbol{\mathcal{L}}_k^{(\alpha \beta)}(t,s)\mathbf{L}_k^{(\beta)^{\texttt{H}}}(t,s)
+\mathbf{L}_k^{(\alpha)}(t,s)  \boldsymbol{\Omega}_k^{(\alpha \beta)}(s)   \mathbf{L}_k^{(\beta)^{\texttt{H}}}(t,s),
\end{multline}
for $ t<s$,  with  initial condition $\boldsymbol{\mathcal{V}}_k^{(\alpha \beta)}(t,t)=\boldsymbol{\mathcal{V}}_k^{(\alpha \beta)}(t)$, computed from \eqref{Kij(t,t)}. 
Moreover,  
\begin{multline}\label{L^ij(t,s) smooth}
\boldsymbol{\mathcal{L}}_k^{(\alpha \beta)}(t,s) =   \left[\mathbf{M}_k^{(\alpha)}(t,s-1)-\mathbf{M}_k^{(\alpha \beta)}(t,s-1)\right]\\ \times \mathbf{A}_k^{\texttt{H}}(s)\mathbf{H}_k^{(\beta)^\texttt{H}}(s),  \quad s>t+1,
\end{multline}
and  $\boldsymbol{\mathcal{L}}_k^{(\alpha \beta)}(t,t+1)=  \mathbf{A}_k^{\texttt{H}}(s)\mathbf{J}_k^{(\alpha\beta)}(t,t+1)$.  $ \mathbf{M}_k^{(\alpha)}(t,s)$ is obtained from \eqref{Mk},  and 
\begin{multline}\label{O^ij(t,s) smooth}
\mathbf{M}_k^{(\alpha \beta)}(t,s) = \mathbf{M}_k^{(\alpha \beta)}(t,s-1)+
\boldsymbol{\mathcal{L}}_k^{(\alpha \beta)}(t,s)\boldsymbol{\Omega}_k^{(\beta)^{-1}}(s)\boldsymbol{{J}}_k^{(\beta)^\texttt{H}}(s)\\ + \mathbf{L}_k^{(\alpha)}(t,s)\boldsymbol{{J}}_k^{(\beta\alpha)^\texttt{H}}(s),  \quad  t<s,
\end{multline}
with  initial condition $\mathbf{M}_k^{(\alpha \beta)}(t,t)=   \mathbf{A}_k(t)\mathbf{Q}^{(\alpha\beta)}(t) $.    
\end{theorem}

\begin{remark}
Note that, for $\alpha=\beta$,   $\boldsymbol{\mathcal{L}}_k^{(\alpha \beta)}(t,s) =\boldsymbol{\mathcal{L}}_k^{(\alpha)}(t,s)=\boldsymbol{0}_{kn\times kn}$.
\end{remark}

\section{Numerical simulations}
This section presents a numerical analysis of the proposed algorithms under the $\mathbb{T}_k$-properness conditions  for $k=1$ and $2$. First, the performance and accuracy of the $\mathbb{T}_k$-proper distributed filtering, prediction, and smoothing algorithms are evaluated within both the $\mathbb{T}_1$ and $\mathbb{T}_2$ settings. Next, their  superior behavior compared to quaternion-domain counterparts is illustrated. Finally, the computational advantages of the $\mathbb{T}_k$-proper algorithms are examined in terms  of the number of observations.

To this end,  consider a tessarine Wiener process ${x}(t)$ whose real vector has a correlation matrix
\begin{equation*}\label{Wiener}
\boldsymbol{\Gamma}_{\mathbf{x}^r}(t,s)=\boldsymbol{\mathcal{W} } \min(t,s), \qquad t,s\geq 0
\end{equation*}
with 
\begin{equation}\label{VarWiener}
\boldsymbol{\mathcal{W}}=\left[\begin{array}{cccc} a_1 & 0 & a_3 & a_4 \\ 0 & a_2 & a_4 & a_3 \\ a_3 & a_4 & a_1 & 0  \\ a_4 & a_3 & 0 & a_2 \end{array}\right].
\end{equation}
Therefore,  $x(t)$ is widely factorizable in the sense of Definition \ref{def2},  with $\bar{\mathbf{A}}(t)=4\boldsymbol{\mathcal{J}}_n \boldsymbol{\mathcal{W} }\boldsymbol{\mathcal{J}}_n^{\texttt{H}}  $ and $\bar{\mathbf{B}}(t)=t \boldsymbol{I}_{4} $.

The signal is assumed to be observed by three sensors in accordance to the equation \eqref{obs},  in which $v^{(\alpha)}(t)=\lambda_{\alpha}u(t)$,  with $\lambda_1=0.2$, $\lambda_2=0.5$, $\lambda_3=0.6$, and $u(t)$ a tessarine Gaussian white noise, with real { pseudo}-autocorrelation matrix
\begin{equation}\label{VarW}
\boldsymbol{\Gamma}_{\mathbf{u}^r}(t,s) =\left[ \begin{array}{cccc} 6 & 0 & 4 & 0 \\ 0 & 6& 0 & 4 \\ 4 & 0 & 6 & 0  \\ 0 & 4 & 0 & 6 \end{array}   \right]\delta_{t,s}.
\end{equation}

With the purpose of analyzing both $\mathbb{T}_1$ and $\mathbb{T}_2$-proper scenarios, the following values are considered in \eqref{VarWiener} and \eqref{VarW}:

\begin{itemize}
\item $\mathbb{T}_1$-proper setting: $a_1=a_2=7.6$, $a_3=-2$,  $a_4=0$.

\item $\mathbb{T}_2$-proper setting: $a_1=5.6$, $a_2=2$,  $a_3=0.6$, $a_4=1.2$.

\end{itemize}

\begin{remark}\label{remEJ}
In the $\mathbb{T}_1$-proper setting,  $x(t)$ is $\mathbb{T}_1$-factorizable,  with $A_1(t)=30.4-8\jmath$ and  $B_1(t)=t$ in \eqref{factk1}.  In a similar way, in the $\mathbb{T}_2$-proper setting,  $x(t)$ is $\mathbb{T}_2$-factorizable, with $A_2(t)=\left[\begin{array}{cc} 15.2+2.4\jmath & 7.2+4.8\jmath \\ 7.2-4.8\jmath & 15.2+2.4\jmath  \end{array}\right]$ and  $B_2(t)=t\boldsymbol{I}_{2}$ in \eqref{factk1}. 
\end{remark}

Moreover,   the following random variables $\gamma_\nu^{(\alpha)}(t)$,  $\nu=r, \imath, \jmath,\kappa$,  are considered to characterize  different types of uncertainty at each sensor: continuous fading measurements  at sensor 1, discrete fading measurements at sensor 2,  and missing measurements at sensor 3.
\begin{itemize}
\item $\mathbb{T}_1$-proper setting: $\gamma_r^{(\alpha)}(t)=\gamma_{ \imath}^{(
\alpha)}(t)=\gamma_{ \jmath}^{(\alpha)}(t)=\gamma_{ \kappa}^{(\alpha)}(t)=\gamma^{(\alpha)}(t)$, $\alpha=1,2,3$, with 
 \begin{description}
\item[\hskip-0.5cm {\bf Sensor 1.}]  $\gamma^{(1)}(t)$ is uniformly distributed  on the interval $[0.2, 0.8]$.
\item[{\hskip-0.5cm  \bf Sensor 2.}] $\gamma^{(2)}(t)$  has probability mass function given by $P[\gamma^{(2)}(t)=0 ]=0.3$,  $P[\gamma^{(2)}(t) =0.5 ]=0.2$,  $P[\gamma^{(2)}(t)=1  ]=0.5$.
\item[{\hskip-0.5cm  \bf Sensor 3.}] 
$\gamma^{(3)}(t)$  is a Bernoulli random variable with parameter $p=0.9$. 
\end{description}
\item $\mathbb{T}_2$-proper setting: $\gamma_r^{(\alpha)}(t)=\gamma_{ \jmath}^{(\alpha)}(t)$,  $\gamma_{ \imath}^{(\alpha)}(t)=\gamma_{ \kappa}^{(\alpha)}(t)$, $\alpha=1,2,3$, with 
\begin{description}
\item[\hskip-0.5cm  {\bf Sensor 1.}]  $\gamma_r^{(1)}(t)$  and  $\gamma_{\imath}^{(1)}(t)$  are uniformly distributed on the intervals $[0.15, 0.45]$  and  $[0.1, 0.7]$, respectively.
\item[\hskip-0.5cm   {\bf Sensor 2.}] $\gamma_r^{(2)}(t)$  and  $\gamma_{\imath}^{(2)}(t)$ has probability mass functions given by $P[\gamma_r^{(2)}(t)=0  ]=0.3$, $P[\gamma_r^{(2)}(t) =0.5 ]=0.2$,  $ P[\gamma_r^{(2)}(t)=1]=0.5$, and $P[\gamma_{\imath}^{(2)}(t)=0 ]=0.1$,  $P[\gamma_{\imath}^{(2)}(t) =0.5 ]=0.6$,   $P[\gamma_{\imath}^{(2)}(t)=1 ]=0.3$.

\item[\hskip-0.5cm   {\bf Sensor 3. }]
$\gamma_r^{(3)}(t)$  and  $\gamma_{\imath}^{(3)}(t)$    are Bernoulli random variables with parameters $p_r=0.8$ and $p_{\imath}=0.7$, respectively. 
\end{description}
\end{itemize}

\begin{figure}[htb] 
\centering
\subfloat[ $\mathbb{T}_1$-proper scenario \label{fig2a}]{
\includegraphics[width=3.5in]{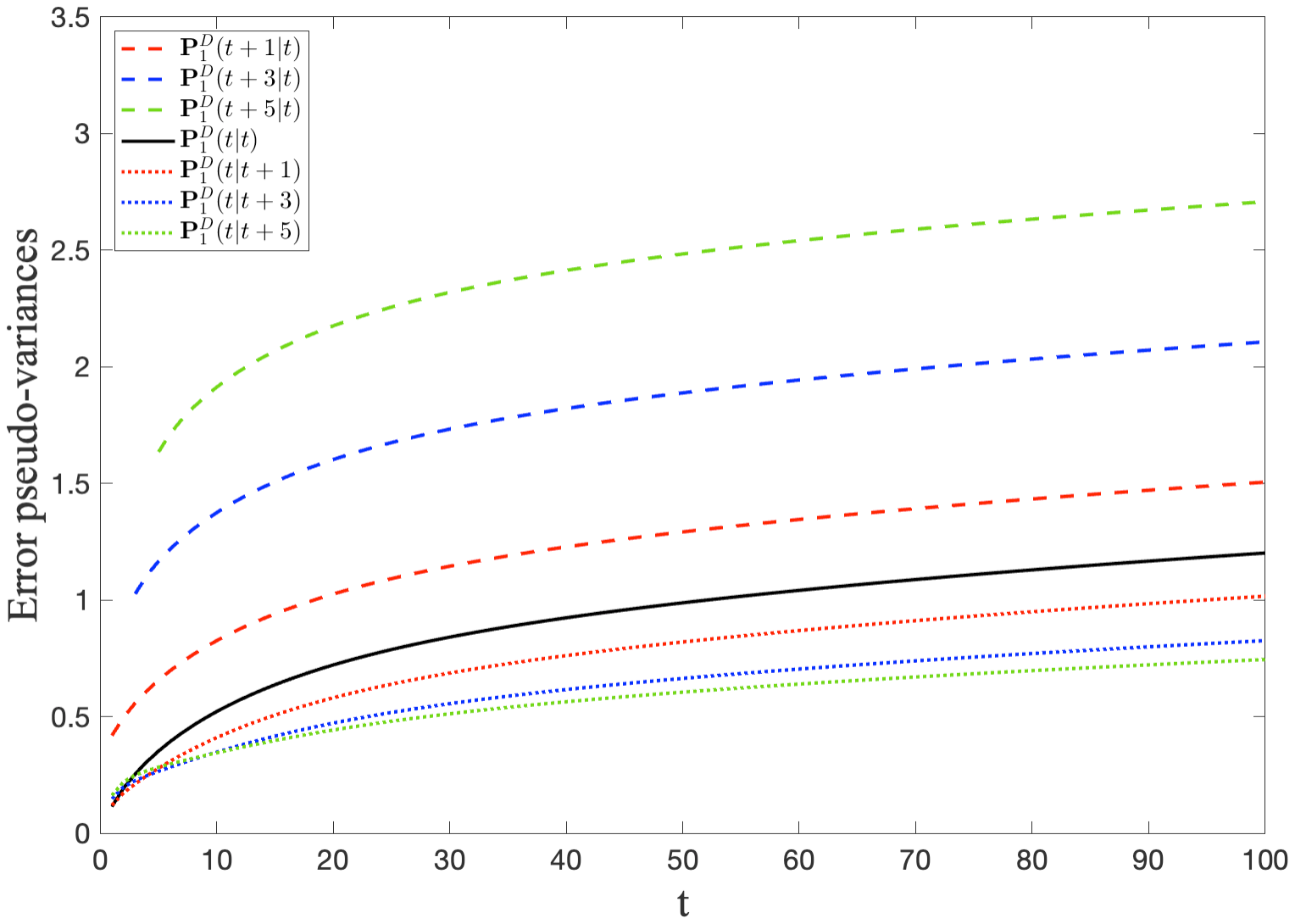} 
}
\\
\subfloat[ $\mathbb{T}_2$-proper scenario \label{fig2b}]{
 \includegraphics[width=3.5in]{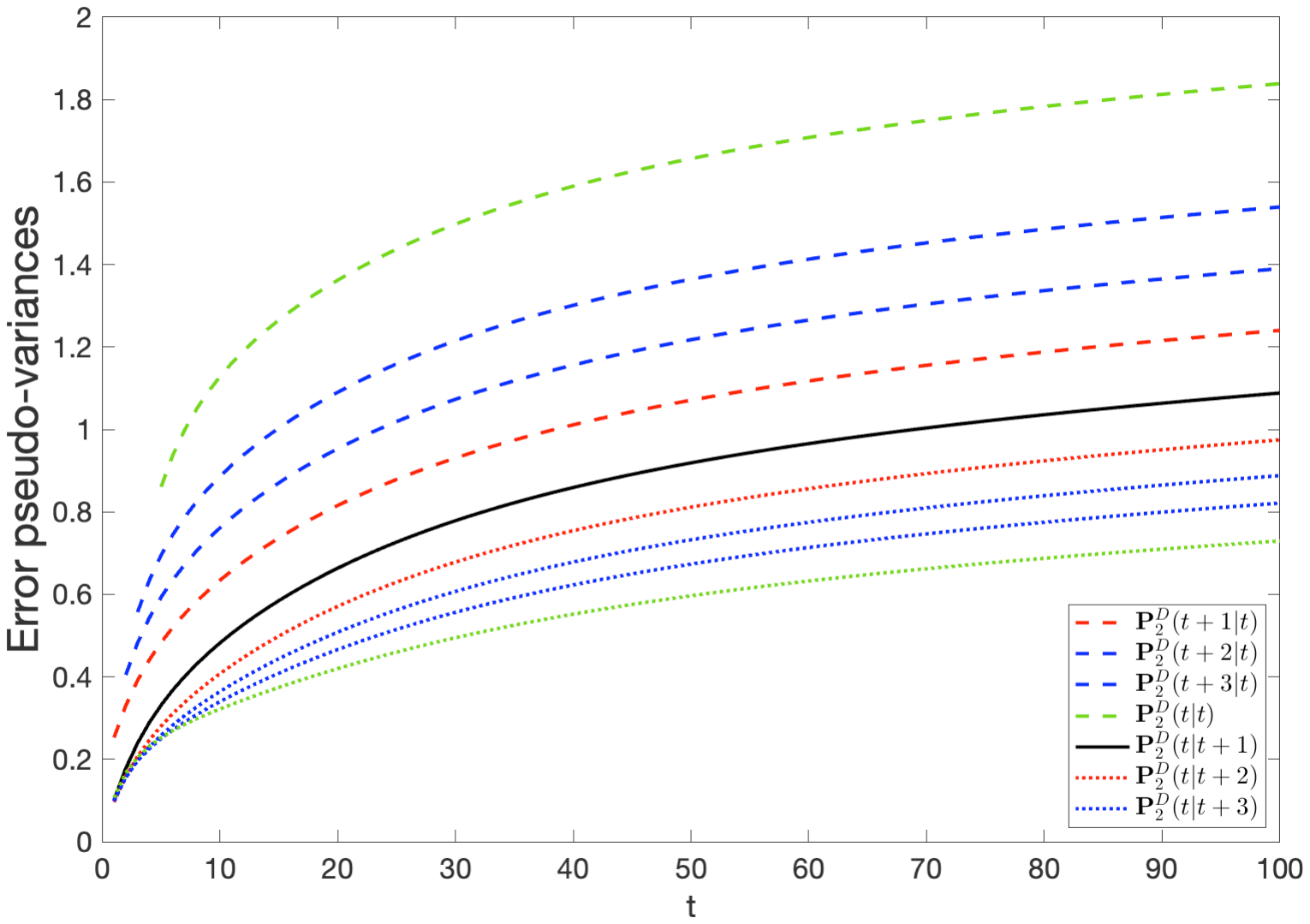}
}
\caption{$\mathbb{T}_k$-proper distributed filtering (solid line), prediction (dashed line) and smoothing (dotted line) error pseudo-variances.}
\end{figure}

To evaluate the behavior of the $\mathbb{T}_k$-proper distributed fusion linear estimation algorithms proposed here, the  error pseudo-variances   ${\mathbf{P}}_k^{D}(t|s)$ have been computed for both  $\mathbb{T}_1$ and  $\mathbb{T}_2$ scenarios. These pseudo-variances are  displayed in Fig.~\ref{fig2a} and Fig.~\ref{fig2b},  respectively, for the finite-step prediction (${\mathbf{P}}_k^{D}(t+\tau|t)$, $\tau=1,3,5$), filtering (${\mathbf{P}}_k^{D}(t|t)$), and finite-step smoothing (${\mathbf{P}}_k^{D}(t|t+\tau)$, $\tau=1,3,5$) problems.  These figures reveal that the prediction error pseudo-variances are greater than those of the  corresponding filters, and they increase as the  prediction time ahead $\tau$ grows. This indicates that,  the performance of the  $\mathbb{T}_k$-proper distributed  prediction estimator is worse than that of the  corresponding filter,  and the performance disparity  grows as the prediction time extends further into the future.  Conversely, the smoothing error pseudo-variances are smaller than those of the corresponding filters, and decrease as the smoothing time lag $\tau$ increases.  Accordingly, the performance of the  $\mathbb{T}_k$-proper distributed  smoothing estimator is  superior to that of the  corresponding filter,  with greater improvements observed as the smoothing time lag $\tau$ increases.

\begin{figure}[htb] 
\centering
\subfloat[ $\mathbb{T}_1$-proper scenario \label{fig1a}]{
 \includegraphics[width=3.5in]{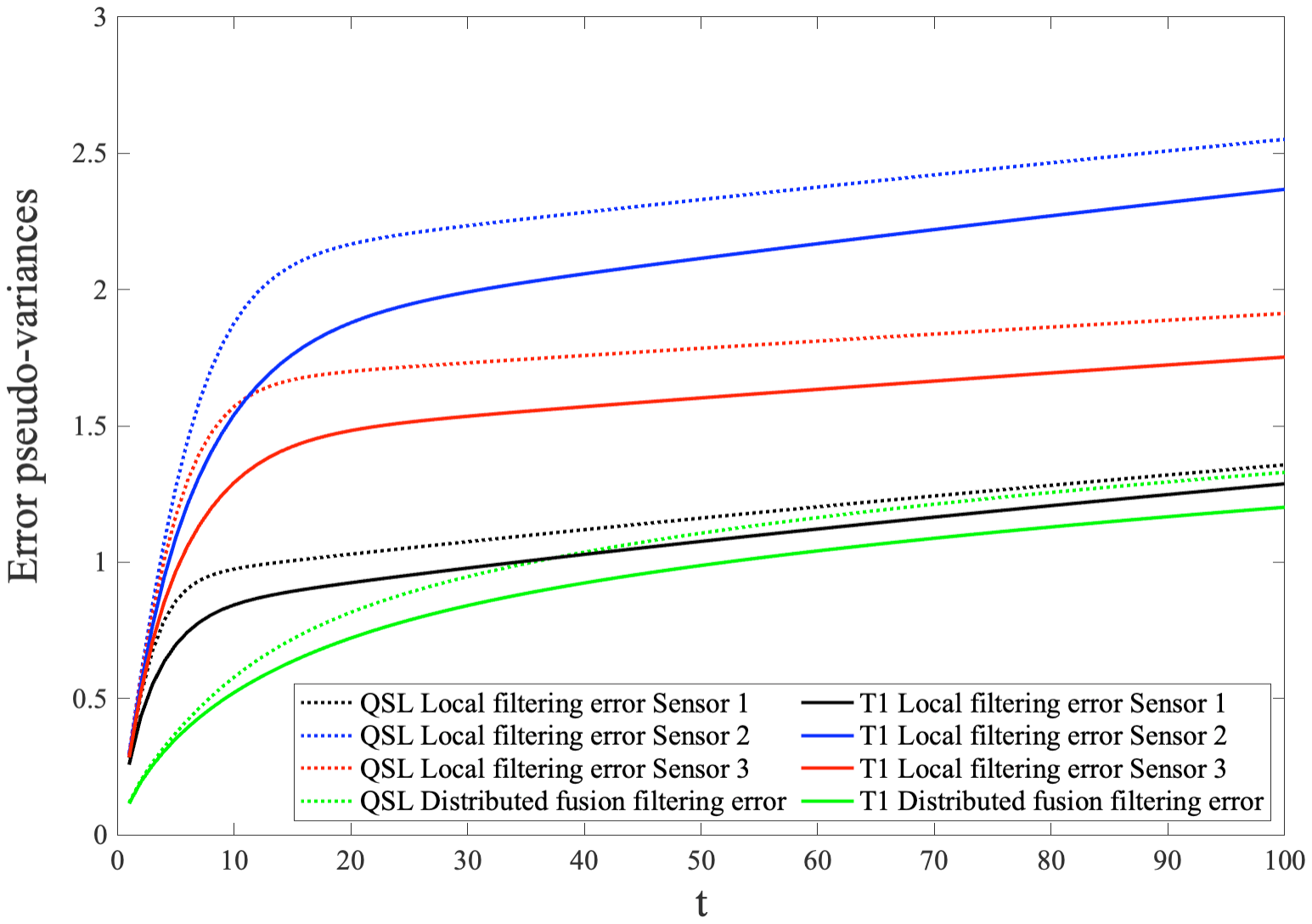} 
}
\\
\subfloat[ $\mathbb{T}_2$-proper scenario \label{fig1b}]{
 \includegraphics[width=3.5in]{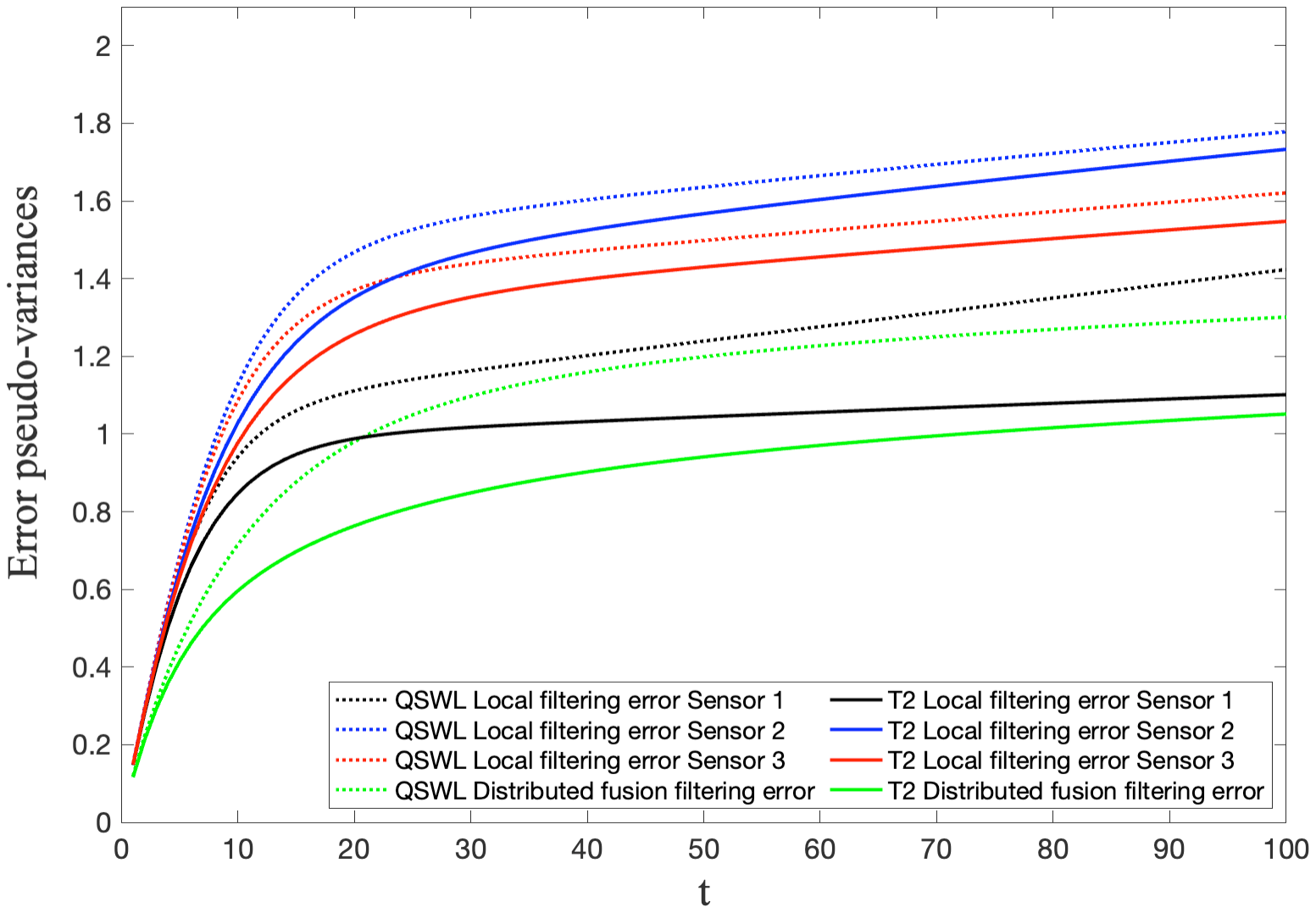}
}
\caption{$\mathbb{T}_1$-proper vs QSL  local and distributed filtering error pseudo-variances in a $\mathbb{T}_1$-proper  setting (a), and $\mathbb{T}_2$-proper vs QSWL  local and distributed filtering error pseudo-variances in a $\mathbb{T}_2$-proper setting (b).}
\end{figure}

The second objective of this analysis is to assess the better performance of the  $\mathbb{T}_k$-proper distributed fusion linear estimation algorithms in comparison to their  quaternion-domain counterparts. To this end, the error pseudo-variances of the $\mathbb{T}_1$-proper and quaternion strictly linear (QSL) local and distributed fusion filters have been computed under $\mathbb{T}_1$-properness conditions. Similarly, the error pseudo-variances  associated with the  $\mathbb{T}_2$-proper and quaternion semi-widely linear (QSWL) local and distributed fusion filters have been evaluated under $\mathbb{T}_2$-properness conditions. 
The outcomes of this comparative study are shown in  Fig. ~\ref{fig1a} and Fig.~\ref{fig1b} for the $\mathbb{T}_1$ and $\mathbb{T}_2$-proper cases, respectively.   These figures illustrate two key findings: Firstly,  the superiority of the distributed filtering estimators over the local filter computed at individual sensors and,  secondly,  the      $\mathbb{T}_k$-proper filtering estimators  consistently outperform their  the quaternion-domain counterparts.

\begin{figure}[htb] 
\centering
\subfloat[ $\mathbb{T}_1$-proper scenario \label{figMET1}]{
 \includegraphics[width=3.5in]{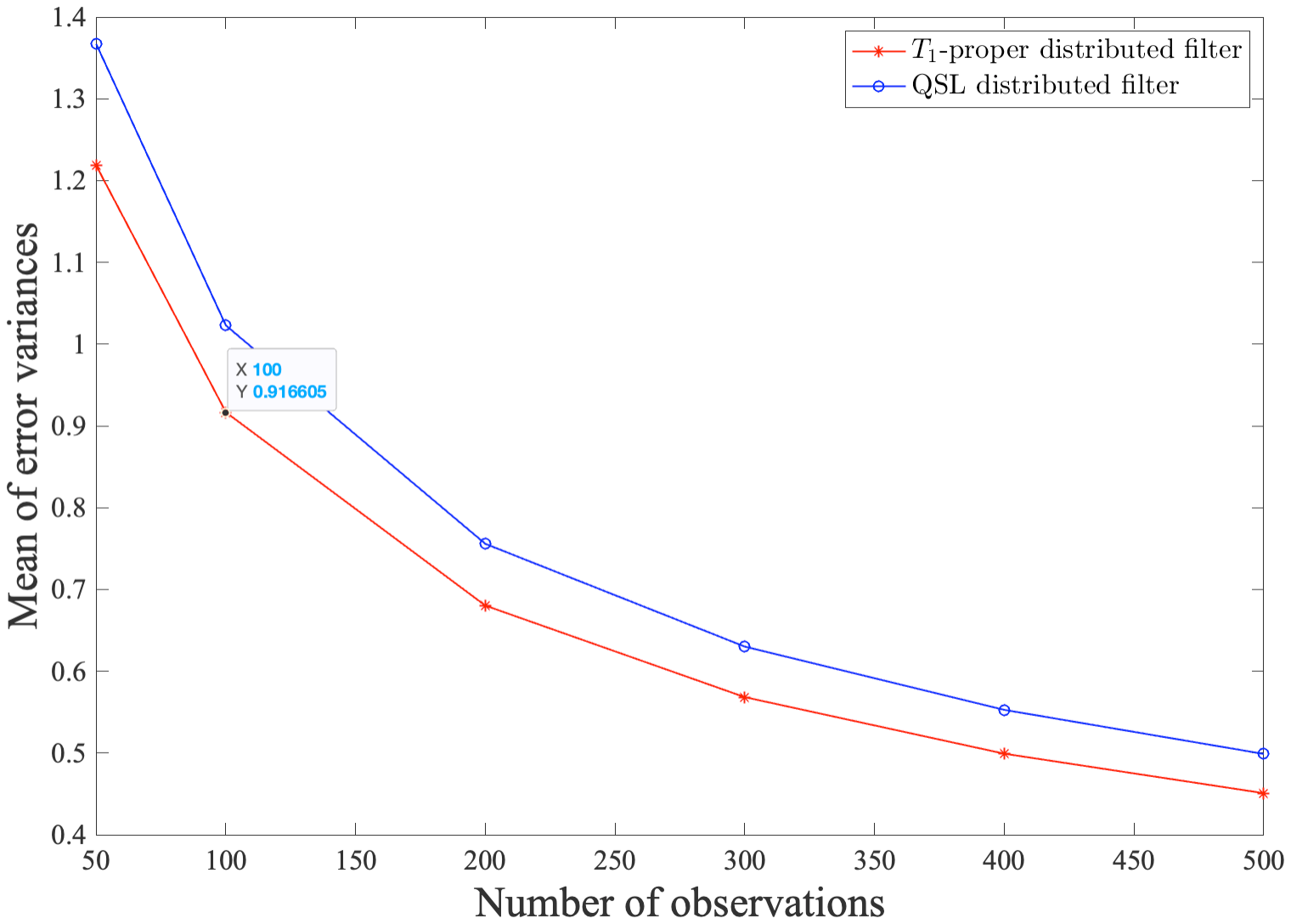} 
}
\\
\subfloat[ $\mathbb{T}_2$-proper scenario \label{figMET2}]{
 \includegraphics[width=3.5in]{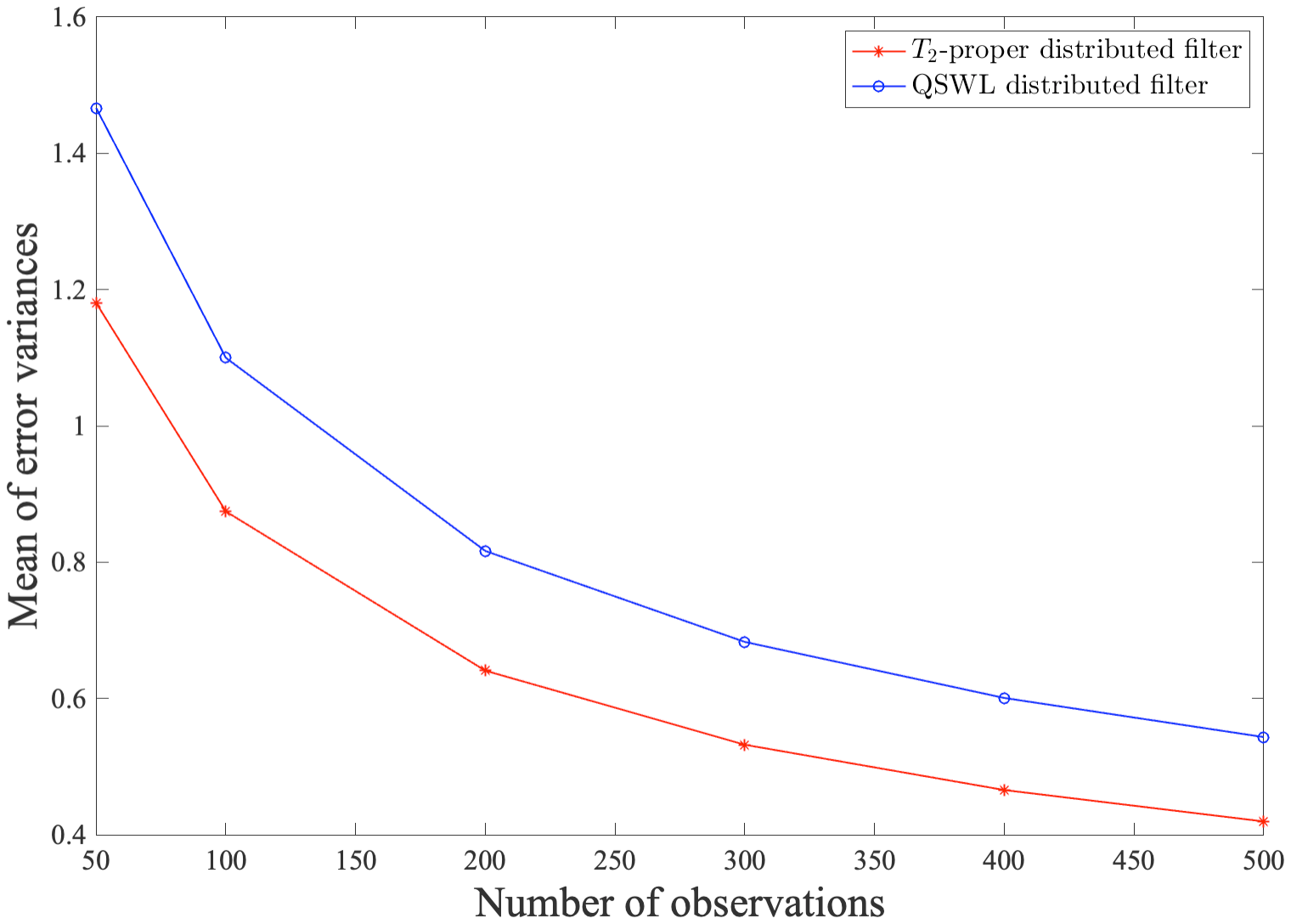} 
}
\caption{Means of the error variances of the  $\mathbb{T}_k$-proper distributed filtering estimates and their counterparts in a $\mathbb{T}_1$-proper setting (a),  and  in a $\mathbb{T}_2$-proper setting (b),  in terms of the number of observations. }
\end{figure}

To evaluate the 
effect of the number of observations on the behavior  of the  $\mathbb{T}_k$-proper distributed linear estimation algorithms compared to their counterparts in the quaternion domain,  the mean of the corresponding distributed fusion filtering error pseudo-variances has been computed  and displayed  in Fig. ~\ref{figMET1} and  Fig. ~\ref{figMET2}, under a $\mathbb{T}_1$ and $\mathbb{T}_2$-proper  scenarios,  respectively.  These figures demonstrate that the accuracy of the distributed fusion filtering estimators increases with the number of observations,
and that the superior performance of the $\mathbb{T}_k$ estimators over their quaternion domain counterparts remains consistent in the steady state.  Similar results are observed for the prediction, and smoothing problems.

\begin{figure}[htb] 
\centering
\subfloat[ $\mathbb{T}_1$-proper scenario \label{fig3a}]{
 \includegraphics[width=3.5in]{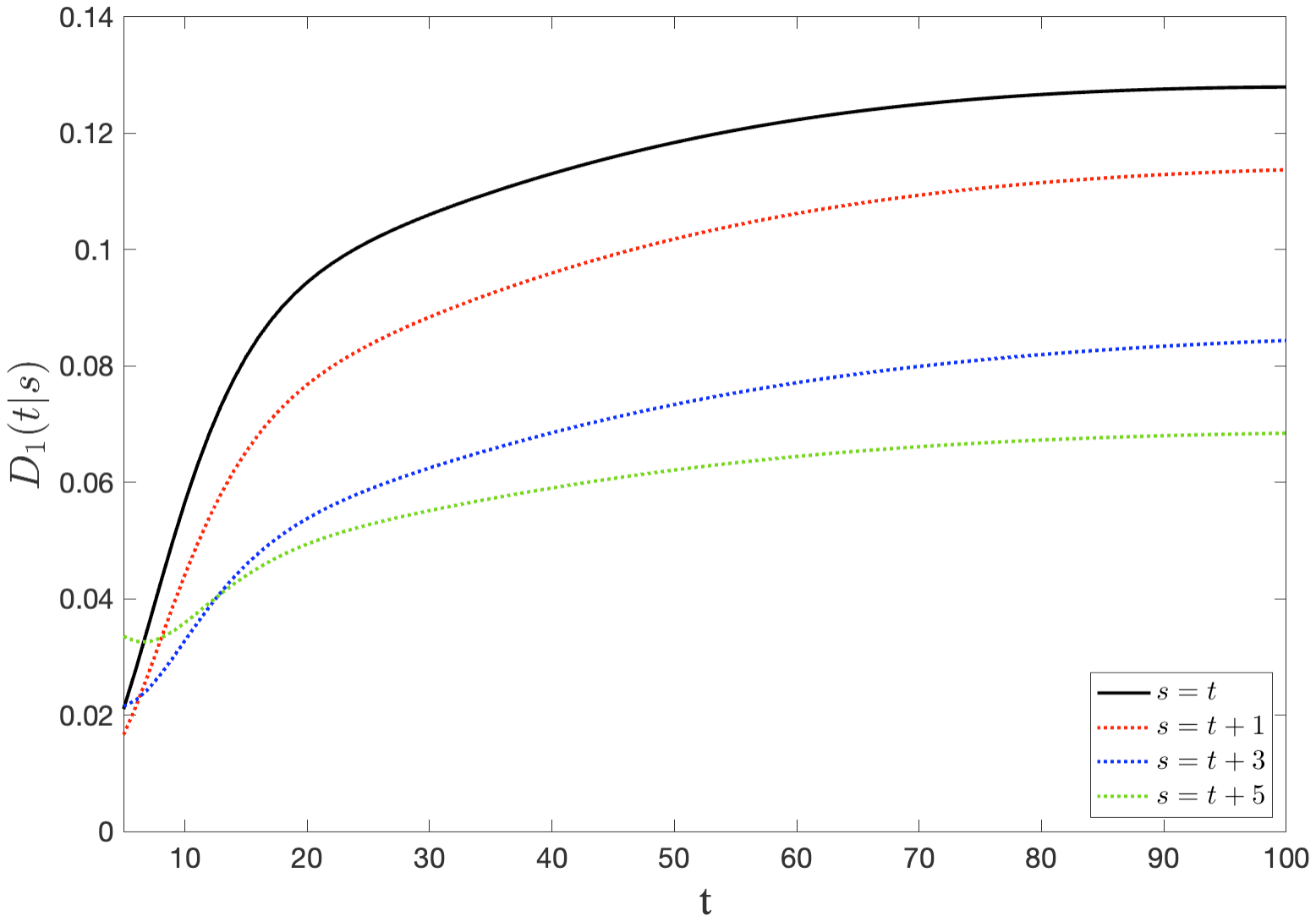} 
}
\\
\subfloat[ $\mathbb{T}_2$-proper scenario \label{fig3b}]{
 \includegraphics[width=3.5in]{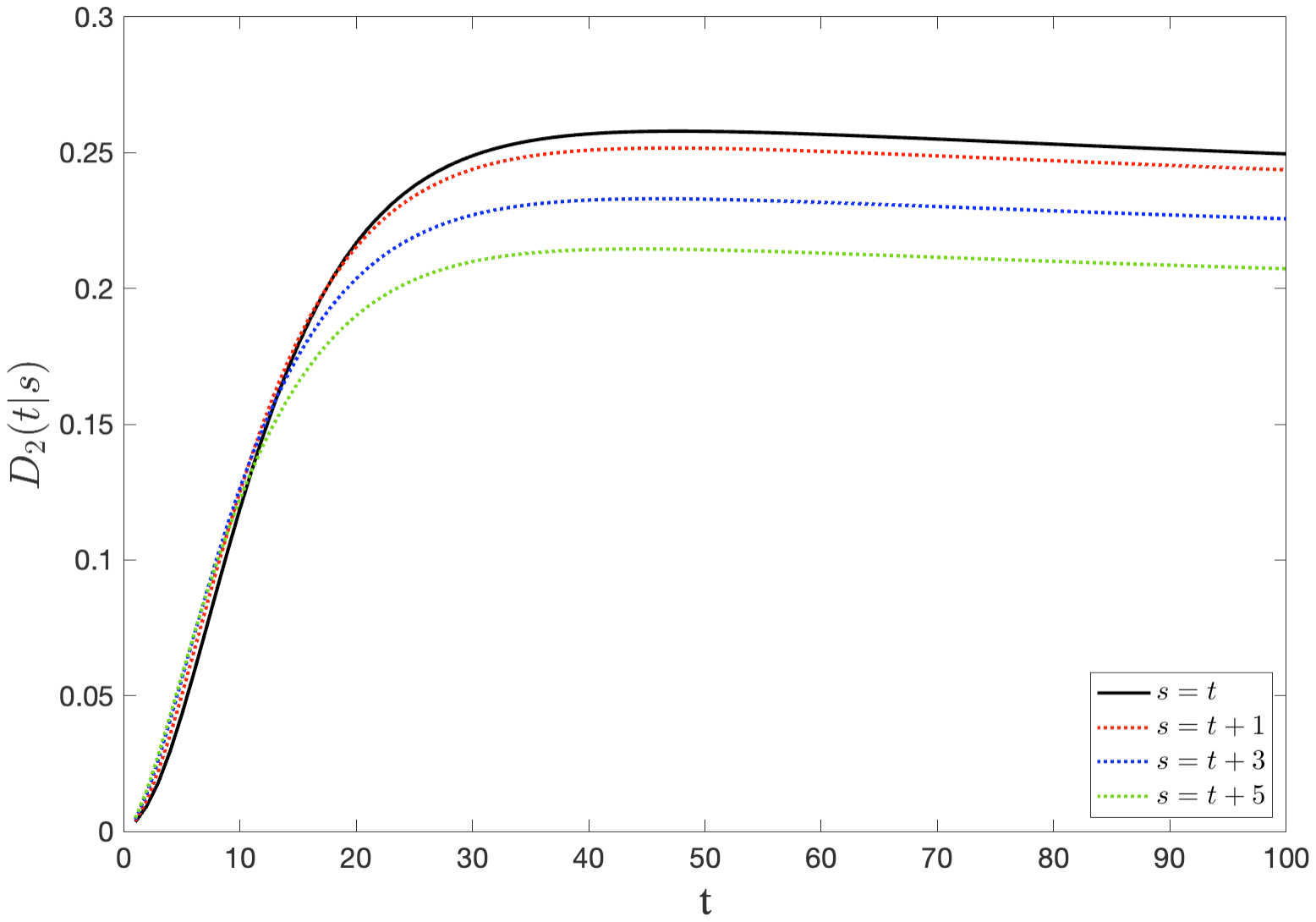} 
}
\caption{Difference between QSL and $\mathbb{T}_1$-proper filtering and smoothing error variances in a $\mathbb{T}_1$-proper setting (a), and between  QSWL  and $\mathbb{T}_2$-proper  filtering and smoothing error variances in a $\mathbb{T}_2$-proper setting }
\end{figure}

Next,  a  comparative analysis is conducted for the distributed fusion prediction and smoothing problems by computing  the following differences:
\begin{align*}
D_1(t|s) & =P_{\texttt{QSL}}^D(t|s)-P_{1}^D(t|s),\\
D_2(t|s) & =P_{\texttt{QSWL}}^D(t|s)-P_{2}^D(t|s),
\end{align*}
where $P_{\texttt{QSL}}^D(t|s)$  and $P_{\texttt{QSWL}}^D(t|s)$ denote the QSL and QSWL distributed error pseudo-variances, respectively. It is worth noting  that these differences are identical for both the filtering and prediction problems. Consequently, only the results for the filtering and smoothing   problems are presented in Fig.~\ref{fig3a} and Fig.~\ref{fig3b}, for the $\mathbb{T}_1$ and $\mathbb{T}_2$-proper cases, respectively.  As  it can be observed, the  $\mathbb{T}_k$-proper distributed fusion estimators  outperform 
their counterparts in the quaternion domain (positive differences),  although this improvement is less evident as the  smoothing time lag increases. 

\begin{figure}[htb] 
\centering
\subfloat[ $\mathbb{T}_1$-proper scenario \label{timeT1}]{
 \includegraphics[width=3.5in]{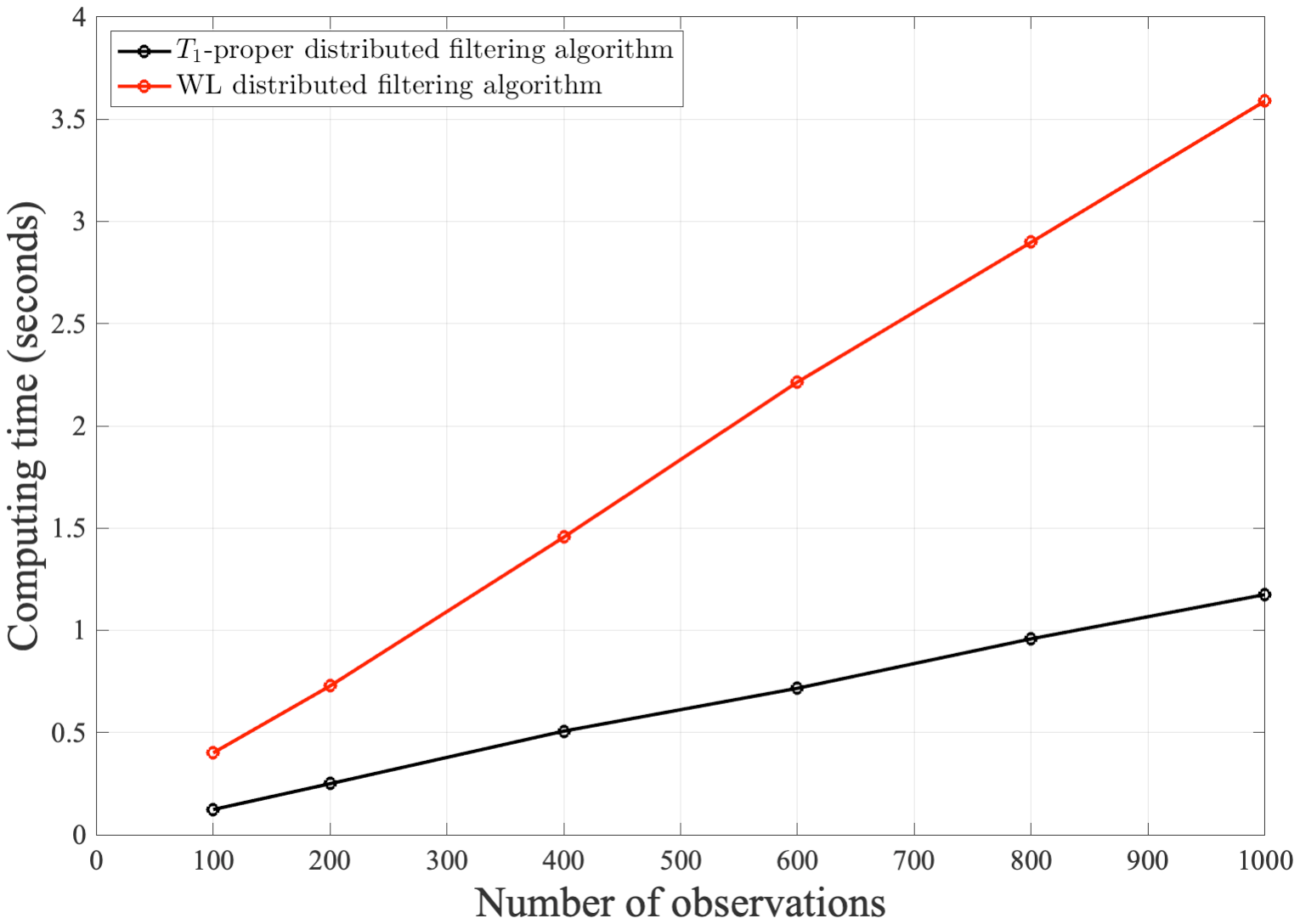} 
}
\hfill
\subfloat[ $\mathbb{T}_2$-proper scenario \label{timeT2}]{
 \includegraphics[width=3.5in]{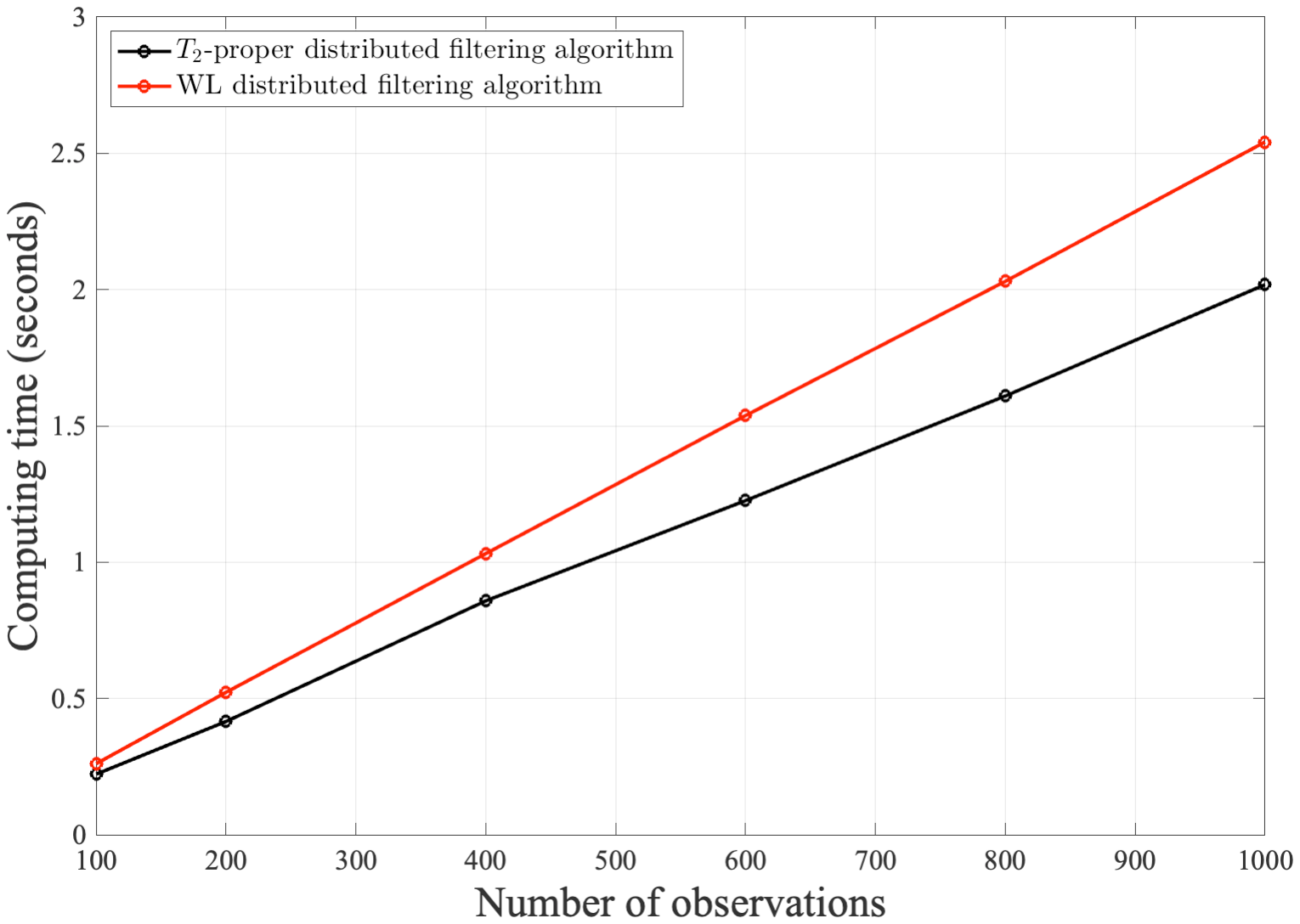} 
}
\caption{Computing time (in seconds) using a $\mathbf{T}_k$-proper and a WL distributed fusion filtering algorithm  in a $\mathbb{T}_1$-proper setting (a),  and  in a $\mathbb{T}_2$-proper setting (b),  and varying the number of observations from $50$ to $500$.}
\end{figure}

In general,  under $\mathbb{T}_k$-properness conditions, the $\mathbb{T}_k$-proper distributed fusion estimators coincide with  those obtained under the optimal (WL) processing, 
whereas a significant computational gain is achieved by reducing the dimensionality to half (in the case of $k=2$) or a quarter (in the case of $k=1$) of the equations used for their calculation.  In Fig.~\ref{timeT1} and Fig.~\ref{timeT2},  the computation time of the $\mathbb{T}_k$-proper distributed filtering   algorithms proposed is compared to that of 
the corresponding  WL distributed filtering  algorithm, in both  $\mathbb{T}_1$- and $\mathbb{T}_2$-proper scenarios, respectively,  as  the number of observations increases.   The computations were performed on an iMac equipped with an Apple M1 chip, which has 8 total cores (4 performance  and 4 efficiency ) and 16 GB of RAM.  The software used for the computations was MATLAB R2021a (9.10.0).   As shown in both figures, the  $\mathbb{T}_k$ -proper algorithms consistently outperform the WL counterparts  in terms of execution time, demonstrating a significant reduction in the computational cost.  Furthermore,   this computational advantage becomes increasingly significant as the number of observations grows, particularly in the $\mathbb{T}_1$-proper case, where the  execution time of the WL algorithm is more than twice that of the $\mathbb{T}_k$-proper  algorithm  for $1,000$ observations.

\section{Conclusions}

This paper presents a novel approach to distributed fusion estimation in multi-sensor networks with fading measurements  for the class of widely factorizable signals.   Specifically,  all types of estimation problems have been addressed: filtering, prediction, and smoothing. 

Unlike most previous results, the problem has been tackled within a hypercomplex (tessarine) framework, under conditions of $\mathbb{T}_k$-properness.  The main benefit of this structure lies in tits ability to reformulate the problem into an equivalent form with a reduced dimensionality -- down to a quarter (for $k=1$) or to half (for $k=2$),  thereby leading to a significant reduction in the computational complexity of the estimators proposed.  

Consequently,   by using  second-order statistical information,  distributed fusion algorithms have been devised to compute filtering, prediction, and smoothing estimators, as well as their mean squared errors. These algorithms are computationally more efficient than the traditional real-field estimators. Moreover, they are valid for stationary and non-stationary signals,  and are applicable in scenarios where a state-space model is not directly accessible,
 as they require only second-order statistical information.  

The behaviour and efficacy  of the $\mathbb{T}_k$-proper estimators  proposed have been verified under different uncertainties scenarios  and  $\mathbb{T}_k$-proper contexts.  The results confirm that $\mathbb{T}_k$-proper estimators outperform their counterparts in the quaternion domain, in a scenario of  $\mathbb{T}_k$-properness.

  It is worth noting that,  in practice,  any 4D hypercomplex algebra  --such as quaternions,  beta quaternions,  and generalized Segre’s quaternions-- could be used,  provided certain proper conditions are met.  In this way,  they enable the best estimators with computational savings not achievable through  real-valued approaches.  Additionally,  statistical tests are available to determine not only whether, under specific criteria, a signal is proper or not, but also to identify the most suitable algebra, if one exists,  within a family that satisfies the proper necessary conditions for dimensionality reduction \cite{Navarro2022, Navarro2024}.
  
\appendices

\section{Proof of Proposition \ref{prop1}}\label{ap-1}
First of all, from \eqref{obs} and Property \ref{propstar},  the following  augmented measurement equation is defined:
\begin{equation}\label{obs4D}
\bar{\mathbf{y}}^{(\alpha)}(t)=\boldsymbol{\mathcal{D}}^{{\boldsymbol{\gamma}}^{(\alpha)}}(t) \bar{\mathbf{x}}(t)+\bar{\mathbf{v}}^{(\alpha)}(t).
\end{equation}
From \eqref{obs4D} and using properties of the Hadamard product, it follows that
\begin{align}\label{covxy}
\boldsymbol{\Gamma}_{\bar{\mathbf{x}}\bar{\mathbf{y}}}^{(\alpha)}(t,s)= & \boldsymbol{\Gamma}_{\bar{\mathbf{x}}}(t,s)\boldsymbol{\Pi}^{(\alpha)}(s),\\ \label{covy}
\boldsymbol{\Gamma}_{\bar{\mathbf{y}}}^{(\alpha)}(t,s)= & \boldsymbol{\Delta}^{(\alpha)}(t,s) +\boldsymbol{\Gamma}_{\bar{\mathbf{v}}}^{(\alpha)}(t,s),
\end{align}
with  $\boldsymbol{\Delta}^{(\alpha)}(t,s)=4\boldsymbol{\mathcal{J}}_n \left\{ \boldsymbol{\Gamma}_{\boldsymbol{\gamma}^r}^{(\alpha)}(t,s)\circ \boldsymbol{\Gamma}_{\mathbf{x}^{r}}(t,s) \right\} \boldsymbol{\mathcal{J}}_n^{\texttt{H}}$,
and $\boldsymbol{\Pi}^{(\alpha)}(s)=\E[\boldsymbol{\mathcal{D}}^{{\boldsymbol{\gamma}}^{(\alpha)}}(s)]$.

In both cases of $T_k$-properness,  $k=1,2$,   the necessary conditions are directly verified by applying the joint $\mathbb{T}_k$-properness definition  and Remark \ref{Rem3} 
on \eqref{covxy}-\eqref{covy}.  

For the sufficient condition, assume that  $\mathbf{x}(t)$ is $\mathbb{T}_1$-proper (respectively, $\mathbb{T}_2$-proper) and  $\mu_{j,\nu}^{(\alpha)}(t)=\mu_{j}^{(\alpha)}(t)$ (respectively, $\mu_{j,\mathrm{r}}^{(\alpha)}(t)=\mu_{j,\jmath}^{(\alpha)}(t)$, $\mu_{j,\imath}^{(\alpha)}(t)=\mu_{j,\kappa}^{(\alpha)}(t)$),  $
j=1,\dots,n$. Then,  from  \eqref{covxy},  it is evident that $\mathbf{x}(t)$ and $\mathbf{y}^{(\alpha)}(t)$ are cross $\mathbb{T}_1$-proper (respectively, $\mathbb{T}_2$-proper).  Furthermore,  since  $\mathbf{v}^{(\alpha)}(t)$ is also $\mathbb{T}_1$-proper (respectively, $\mathbb{T}_2$-proper),  and $\sigma_{j,\nu}^{(\alpha)^2}(t)=\sigma_{j}^{(\alpha)^2}(t)$ (respectively,  $\sigma_{j,\mathrm{r}}^{(\alpha)^2}(t)=\sigma_{j,\jmath}^{(\alpha)^2}(t)$, $\sigma_{j,\imath}^{(\alpha)^2}(t)=\sigma_{j,\kappa}^{(\alpha)^2}(t)$),  $
 j=1,\dots,n$,  it is straightforward to verify, from \eqref{covy}, that $\mathbf{y}^{(\alpha)}(t)$ is also $\mathbb{T}_1$-proper (respectively, $\mathbb{T}_2$-proper).

\section{Proof of Proposition \ref{prop2}}\label{ap-2}
Denoting  $\boldsymbol{\mathcal{J}}_k=\mathcal{C}_k\boldsymbol{\mathcal{J}}_n$,  $\mathbf{x}_k(t)$ and $\mathbf{y}_k^{(\alpha)}(t)$ can be expressed as  
\begin{align}\label{eqxk}
\mathbf{x}_k(t)= & 2\boldsymbol{\mathcal{J}}_k \mathbf{x}^r(t),\\
\mathbf{y}_k^{(\alpha)}(t)= & 2\boldsymbol{\mathcal{J}}_k \diag\left(\boldsymbol{\gamma}^{(\alpha)^r}(t)  \right) \mathbf{x}^r(t)+\mathbf{v}_k^{(\alpha)}(t).
\label{eqyk}
\end{align}
Moreover, from Proposition \ref{prop1}, it follows that 
$\boldsymbol{\mathcal{J}}_k \E\left[\diag\left(\boldsymbol{\gamma}^{(\alpha)^r}(t)  \right)\right]=\mathbf{H}_k^{(\alpha)}(t)\boldsymbol{\mathcal{J}}_k $. Consequently, it is not difficult to check that $\boldsymbol{\Gamma}_{\mathbf{y}_k\mathbf{x}_k}^{(\alpha)}(t,s)=\boldsymbol{\Gamma}_{\mathbf{z}_k\mathbf{x}_k}^{(\alpha)}(t,s)$ and, applying  Hadamard product properties,  $\boldsymbol{\Gamma}_{\mathbf{y}_k}^{(\alpha\beta)}(t,s)=\boldsymbol{\Gamma}_{\mathbf{z}_k}^{(\alpha)}(t,s)$.

\section{Proof of Theorem \ref{ThLocalFilter}}\label{ap1}

At each sensor $\alpha=1,\dots,R$,   the local linear MMSE estimator $\hat{\mathbf{x}}_k^{(\alpha)}(t|t)$ of  $\mathbf{x}_k(t)$  is defined as the projection of $\mathbf{x}_k(t)$  onto the set of observations $\{ \mathbf{y}_k^{(\alpha)}(1),\ldots, \mathbf{y}_k^{(\alpha)}(t)\}$.   Note that, although the set of tessarine random variables lacks a Hilbert space structure,  recent advancements have established  a metric space framework within the tessarine domain that guarantees the existence and uniqueness of the orthogonal projection \cite{Navarro2021}.

The proof is based on an innovation approach.  
Then,  for any $t,s\geq 1$,  
the local linear MMSE estimator  
$\hat{\mathbf{x}}_k^{(\alpha)}(t|s)$ 
 is obtained by  projecting $\mathbf{x}_k(t)$ onto the innovations $\left\{{\boldsymbol{\varepsilon}}_k^{(\alpha)}(1),\ldots, {\boldsymbol{\varepsilon}}_k^{(\alpha)}(s)\right\}$,  with ${\boldsymbol{\varepsilon}}_k^{(\alpha)}(s)=\mathbf{y}_k^{(\alpha)}(s)-\hat{\mathbf{y}}_k^{(\alpha)}(s|s-1)$, and $\hat{\mathbf{y}}_k^{(\alpha)}(s|s-1)$  the local linear MMSE one-stage predictor of $\mathbf{y}_k^{(\alpha)}(s)$.  This estimator  
can be expressed as 
\begin{equation}\label{estfilter}
\hat{\mathbf{x}}_k^{(\alpha)}(t|s)=\sum_{j=1}^s   \mathbf{L}_k^{(\alpha)}(t,j){\boldsymbol{\varepsilon}}_k^{(\alpha)}(j), \quad t, s \geq 1
\end{equation}
where $ \mathbf{L}_k^{(\alpha)}(t,j)={\boldsymbol{\Phi}}_k^{(\alpha)}(t,j){\boldsymbol{\Omega}}_k^{(\alpha)^{-1}}(j)$, with ${\boldsymbol{\Phi}}_k^{(\alpha)}(t,j)= \boldsymbol{\Gamma}_{\mathbf{x}_k\boldsymbol{\varepsilon_k}}^{(\alpha)}(t,j)$,
and ${\boldsymbol{\Omega}}_k^{(\alpha)}(j)=\boldsymbol{\Gamma}_{\boldsymbol{\varepsilon_k}}^{(\alpha)}(j,j)$.

In a similar way,  and using \eqref{eqxk},  \eqref{eqyk},  and Proposition  \ref{prop1}, it follows that
\begin{multline}\label{predy}
\hat{\mathbf{y}}_k^{(\alpha)}(t|t-1)=\sum_{j=1}^{t-1}  
\boldsymbol{\Gamma}_{\mathbf{y}_k\boldsymbol{\varepsilon_k}}^{(\alpha)}(t,j){\boldsymbol{\Omega}}_k^{(\alpha)^{-1}}(j) {\boldsymbol{\varepsilon}}_k^{(\alpha)}(j)\\=\mathbf{H}_k^{(\alpha)}(t)\hat{\mathbf{x}}_k^{(\alpha)}(t|t-1).
\end{multline}

Now,  taking $s=t-1$ in \eqref{estfilter}, and using \eqref{Innok},  \eqref{eqxk}, and Proposition  \ref{prop1},  it follows that 
\begin{equation*}\label{Phi}
\begin{split}
{\boldsymbol{\Phi}}_k^{(\alpha)}(t,j){=} & \boldsymbol{\Gamma}_{\mathbf{x}_k\mathbf{z}_k}^{(\alpha)}(t,j) \\ & - \sum_{l=1}^{j-1}{\boldsymbol{\Phi}}_k^{(\alpha)}(t,l)  {\boldsymbol{\Omega}}_k^{(\alpha)^{-1}}(l) {\boldsymbol{\Phi}}_k^{(\alpha)^{\texttt{H}}}(j,l) \mathbf{H}_k^{(\alpha)^{\texttt{H}}}(j)\\
\underset{\eqref{factk1}, \eqref{zk}}{=} & \mathbf{A}_k(t)\mathbf{B}_k^{\texttt{H}}(j)\mathbf{H}_k^{(\alpha)^{\texttt{H}}}(j)\\  &  - \sum_{l=1}^{j-1}{\boldsymbol{\Phi}}_k^{(\alpha)}(t,l)  {\boldsymbol{\Omega}}_k^{(\alpha)^{-1}}(l) {\boldsymbol{\Phi}}_k^{(\alpha)^{\texttt{H}}}(j,l) \mathbf{H}_k^{(\alpha)^{\texttt{H}}}(j).
\end{split}
\end{equation*}

Introducing the auxiliary function 
\begin{multline}\label{eqJ}
\mathbf{J}_k^{(\alpha)}(t)= \mathbf{B}_k^{\texttt{H}}(t)\mathbf{H}_k^{(\alpha)^{\texttt{H}}}(t)\\ - \sum_{l=1}^{t-1}\mathbf{J}_k^{(\alpha)}(l) {\boldsymbol{\Omega}}_k^{(\alpha)^{-1}}(l) {\boldsymbol{\Phi}}_k^{(\alpha)^{\texttt{H}}}(t,l) \mathbf{H}_k^{(\alpha)^{\texttt{H}}}(t), \quad t\geq 2
\end{multline}
with $\mathbf{J}_k^{(\alpha)}(1)=  \mathbf{B}_k^{\texttt{H}}(1)\mathbf{H}_k^{(\alpha)^{\texttt{H}}}(1) $,  it follows that
\begin{equation}\label{eqPhi}
{\boldsymbol{\Phi}}_k^{(\alpha)}(t,j)=\mathbf{A}_k(t) \mathbf{J}_k^{(\alpha)}(j), \quad  1\leq j\leq t.
\end{equation}
Then, defining 
\begin{equation}\label{defE}
\mathbf{e}^{(\alpha)}(t)= \sum_{j=1}^{t}\mathbf{J}_k^{(\alpha)}(j) {\boldsymbol{\Omega}}_k^{(\alpha)^{-1}}(j) {\boldsymbol{\varepsilon}}_k^{(\alpha)}(j), \quad t\geq 1,
\end{equation}
the  expression \eqref{local filter}  for the filter $\hat{\mathbf{x}}_k^{(\alpha)}(t|t)$  holds,  where \eqref{ek} is a direct consequence of 
 \eqref{defE}. 
Similarly,  it is verified that  $\hat{\mathbf{x}}_k^{(\alpha)}(t|t-1)=\mathbf{A}_k(t)\mathbf{e}^{(\alpha)}(t-1)$,  which lead to the expression  \eqref{Innok} for the innovation ${\boldsymbol{\varepsilon}}_k^{(\alpha)}(t)$.

In addition,  the expression \eqref{Jk} for $\mathbf{J}_k^{(\alpha)}(t)$ is obtained by substituting  \eqref{eqPhi} in \eqref{eqJ},  and introducing the matrix
\begin{equation}\label{EqQ}
\mathbf{Q}^{(\alpha)}(t)= \boldsymbol{\Gamma}_{\mathbf{e}}^{(\alpha)} (t,t) \underset{\eqref{defE}}{=}\sum_{j=1}^{t}\mathbf{J}_k^{(\alpha)}(j) {\boldsymbol{\Omega}}_k^{(\alpha)^{-1}}(j) \mathbf{J}_k^{(\alpha)^{\texttt{H}}}(j),
\end{equation}
for $t\geq 2$.  Furthermore,  from \eqref{EqQ},  the formula  \eqref{Qk} is directly devised.

Next,  from Proposition \ref{prop2}, \eqref{zk}, \eqref{predy}, \eqref{EqQ} the pseudo-variance matrix of the innovations $ {\boldsymbol{\varepsilon}}_k^{(\alpha)}(t)$ takes the form 
\begin{multline*}\label{EqO}
\boldsymbol{\Omega}_k^{(\alpha)}(t){=} \mathbf{H}_k^{(\alpha)}(t)  \boldsymbol{\Gamma}_{\mathbf{x}_k}^{(\alpha)} (t,t)\mathbf{H}_k^{(\alpha)^{\texttt{H}}}(t) + \mathbf{R}_k^{(\alpha)}(t)+\boldsymbol{\Sigma}_k^{(\alpha)}(t) \\
 - \mathbf{H}_k^{(\alpha)}(t) \mathbf{A}_k(t) \mathbf{Q}^{(\alpha)}(t-1) \mathbf{A}_k^{\texttt{H}}(t) \mathbf{H}_k^{(\alpha)^{\texttt{H}}}(t).
\end{multline*}
Then,  \eqref{Omega} follows from \eqref{factk1} and \eqref{Jk}.

Finally,  from \eqref{factk1},  \eqref{local filter},  and \eqref{EqQ},   it is easy to check that the pseudo-variance matrix of the filtering errors  $\boldsymbol{\epsilon}_k^{(\alpha)}(t)=\mathbf{x}_k(t)-\hat{\mathbf{x}}_k^{(\alpha)}(t|t)$,  $\mathbf{P}_k^{(\alpha)}(t|t)= \boldsymbol{\Gamma}_{\boldsymbol{\epsilon}_k}^{(\alpha)}(t,t) 
$,  satisfies equation \eqref{P(t|t) local filter}.

\section{Proof of Theorem \ref{LocalPredictor}}\label{ap2}
By similar reasoning as in Theorem 1, the result follows.

\section{Proof of Theorem \ref{LocalSmoother}}\label{ap3}

From \eqref{estfilter}, it is clear that
\begin{equation*}\label{f-i smoother proof}
\hat{\mathbf{x}}_k^{(\alpha)}(t|s)=
\hat{\mathbf{x}}_k^{(\alpha)}(t|s-1)+\mathbf{L}_k^{(\alpha)}(t,s){\boldsymbol{\varepsilon}}_k^{(\alpha)}(j), \quad t < s
 \end{equation*}
 where, following a similar reasoning that in the Proof of Theorem 1, 
 \begin{equation*}\label{Phi2}
{\boldsymbol{\Phi}}_k^{(\alpha)}(t,s)=\mathbf{B}_k(t)\mathbf{A}_k^{\texttt{H}}(s)\mathbf{H}_k^{(\alpha)^{\texttt{H}}}(s)- \boldsymbol{\Gamma}_{\mathbf{x}_k\mathbf{e}}^{(\alpha)}(t,s-1)  \mathbf{H}_k^{(\alpha)^{\texttt{H}}}(s) 
\end{equation*}
and then, denoting $ \mathbf{M}_k^{(\alpha)}(t,s) =\boldsymbol{\Gamma}_{\mathbf{x}_k\mathbf{e}}^{(\alpha)}(t,s) $, equation \eqref{Lk} holds.

Moreover, from \eqref{ek}, the equation \eqref{Mk} for $ \mathbf{M}_k^{(\alpha)}(t,s)$ is easily devised, with the initialization at s=t given by the expression
\begin{multline*}
\mathbf{M}_k^{(\alpha)}(t,t)\underset{\eqref{defE}}=  \sum_{j=1}^t{\boldsymbol{\Phi}}_k^{(\alpha)}(t,j) {\boldsymbol{\Omega}}_k^{(\alpha)^{-1}}(j) \mathbf{J}_k^{(\alpha)^{\texttt{H}}}(j)\\ \underset{\eqref{Qk}, \eqref{eqPhi}}=   \mathbf{A}_k(t)\mathbf{Q}^{(\alpha)}(t).
\end{multline*}

 Finally, equation \eqref{P(t|N) local smoother}  for the smoothing error pseudo-variance matrix of $\hat{\mathbf{x}}_k^{(\alpha)}(t|s)$  is directly derived from \eqref{local smoother}.

\section{Proof of Theorem \ref{teo distributed filter}}\label{proofteo4}
Introduce  the functions $\mathbf{J}_k^{(\alpha \beta)}(t,s) =  \boldsymbol{\Gamma}_{\mathbf{e} \boldsymbol{\varepsilon}_k}^{(\alpha \beta)}(t,s) $, $\boldsymbol{\Omega}_k^{(\alpha \beta)}(t,s) =  \boldsymbol{\Gamma}_{ \boldsymbol{\varepsilon}_k}^{(\alpha \beta)}(t,s) $, $\mathbf{Q}^{(\alpha \beta)}(t,s) =  \boldsymbol{\Gamma}_{\mathbf{e}}^{(\alpha \beta)}(t,s) $. 
Then,  the expression \eqref{Kij(t,t)} is directly obtained  from \eqref{local filter}, with $\mathbf{Q}^{(\alpha \beta)}(t)=\mathbf{Q}^{(\alpha \beta)}(t,t) $. 

Now, since $\boldsymbol{\mathcal{J}}_k \E\left[\diag\left(\boldsymbol{\gamma}^{(\alpha)^r}(t)  \right)\right]=\mathbf{H}_k^{(\alpha)}(t)\boldsymbol{\mathcal{J}}_k $,  it follows that
\begin{equation*}\label{Gamma3}
\boldsymbol{\Gamma}_{\mathbf{y}_k\mathbf{e} }^{(\alpha \beta)}(t,s)=\boldsymbol{\Gamma}_{\mathbf{z}_k\mathbf{e} }^{(\alpha \beta)}(t,s).
\end{equation*}
As a consequence, from \eqref{Innok} and \eqref{zk},  
 it can be deduced that
\begin{multline*}
\mathbf{J}_k^{(\alpha \beta)}(t-1,t)= \boldsymbol{\Gamma}_{\mathbf{e} \mathbf{x}_k}^{(\alpha)}(t-1,t) \mathbf{H}_k^{(\beta)^{\texttt{H}}}(t)\\ -\mathbf{Q}^{(\alpha \beta)}(t-1) \mathbf{A}_k^{\texttt{H}}(t)\mathbf{H}_k^{(\beta)^{\texttt{H}}}(t).
\end{multline*}
Hence,  the formula \eqref{J^ij(t-1,t)} is devised by considering that $\mathbf{x}_k(t)-\hat{\mathbf{x}}^{(\alpha)}_k(t|t-1)\perp \boldsymbol{\varepsilon}_k^{(\alpha)}(j)$,  for $j\leq t-1$,  and the expression \eqref{local predictor} for $s=t-1$.  
Additionally,  denoting $\boldsymbol{\Omega}_k^{(\alpha \beta)}(t)=\boldsymbol{\Omega}_k^{(\alpha \beta)}(t, t)$,  it can be checked that
\begin{multline*}\label{Oij_proof}
\boldsymbol{\Omega}_k^{(\alpha \beta)}(t)=\boldsymbol{\Gamma}_{\mathbf{y}_k}^{(\alpha\beta)}(t,t)\\ - \mathbf{H}_k^{(\alpha)}(t)\mathbf{A}_k (t)\mathbf{Q}^{(\beta)}(t-1) \mathbf{A}_k^{\texttt{H}} (t)\mathbf{H}_k^{(\beta)^\texttt{H}}(t)\\ -  \mathbf{H}_k^{(\alpha)}(t)\mathbf{A}_k (t)\mathbf{J}_k^{(\alpha \beta)}(t-1,t)
\end{multline*}
and thus, from Proposition \ref{prop2},  \eqref{factk1}, \eqref{zk}, and \eqref{J^ij(t)}, the equation \eqref{Omega^ij(t)} holds.

Moreover,  from \eqref{ek}, it is easy to check that $\mathbf{J}_k^{(\alpha \beta)}(t)=\mathbf{J}_k^{(\alpha \beta)}(t,t)$ satisfies the expression \eqref{J^ij(t)}.  And also
\begin{multline*}
\mathbf{Q}^{(\alpha\beta)}(t) =  \mathbf{Q}^{(\alpha\beta)}(t-1)\\ 
+ \mathbf{J}_k^{(\alpha)}(t)\boldsymbol{\Omega}_k^{(\alpha)^{-1}}(t)\boldsymbol{\Omega}_k^{(\alpha\beta)}(t)\boldsymbol{\Omega}_k^{(\beta)^{-1}}(t)\mathbf{J}_k^{(\beta)^{\texttt{H}}}(t)\\ 
+ \mathbf{J}_k^{(\alpha\beta)}(t-1,t)\boldsymbol{\Omega}_k^{(\beta)^{-1}}(t)\mathbf{J}_k^{(\beta)^{\texttt{H}}}(t)\\
+ \mathbf{J}_k^{(\alpha)}(t)\boldsymbol{\Omega}_k^{(\alpha)^{-1}}(t) \mathbf{J}_k^{(\beta\alpha)^{\texttt{H}}}(t-1,t)
\end{multline*}
 which leads to the expression \eqref{Qij(t)}, by virtue of \eqref{J^ij(t)}.

\section{Proof of Theorem \ref{teo distributed predictor}}\label{proofteo5}
The proof follows from \eqref{local predictor}.

\section{Proof of Theorem \ref{teo distributed smoother}}\label{proofteo6}
Firstly,  introducing the function 
\begin{equation}\label{Lij}
\boldsymbol{\mathcal{L}}_k^{(\alpha \beta)}(t,s)=E\left[ \hat{\mathbf{x}}_k^{(\alpha)}(t|s-1)\boldsymbol{\varepsilon}_k^{(\beta)^\texttt{H}}(s)  \right],
\end{equation}
the expression \eqref{K_k^(ij) smooth}  is directly devised from \eqref{local smoother}.

Moreover,   the expression \eqref{L^ij(t,s) smooth}  for $\boldsymbol{\mathcal{L}}_k^{(\alpha \beta)}(t,s)$ is obtained by using  \eqref{Innok}  on \eqref{Lij}, and introducing the function
\begin{equation}\label{OOij}
\mathbf{M}_k^{(\alpha \beta)}(t,s)=E\left[ \hat{\mathbf{x}}_k^{(\alpha)}(t|s)\mathbf{e}^{(\beta)^\texttt{H}}(s)  \right].
\end{equation}

Now,  using \eqref{ek} and \eqref{local smoother} in \eqref{OOij},  the recursive equation \eqref{O^ij(t,s) smooth}  holds.

%

%
%
%
%
%
%


\nocite{*}
\bibliographystyle{IEEE}

\begin{thebibliography}{00}\label{references}
\bibitem{Singh2024}
A.K.  Singh, S. Berretti,   
{\it Data Fusion Techniques and Applications for Smart Healthcare, } 
1$^{st}$ ed.  Academic Press,  2024.   

\bibitem{Kurkin2017}
A.A.  Kurkin,  D.Y.  Tyugin,  V.D.  Kuzin,  A.G.  Chernov,  V.S.  Makarov,  P.O.  Beresnev,  V.I.  Filatov,   D.V.  Zeziulin,     
``Autonomous mobile robotic system for environment monitoring in a coastal zone, ''
{\it  Procedia Comput.   Sci.},  {vol. 103},  pp.  459--465,  {Dec. 2017},   doi:10.1016/j.procs.2017.01.022. 


%


\bibitem{Zewge2023}
{N.S.  Zewge,  H.A.  Bang,  
``A distributionally robust fusion framework for autonomous multisensor spacecraft navigation during entry phase of mars entry, descent, and landing,''  
 {\it Remote Sens.},  {vol.15},  no.4,  p.1139,  Feb.  2023,  doi:10.3390/rs15041139.
}  




\bibitem{Yan 2019}
J.  Yan, Z.  Xu, X.  Luo, C.  Chen, and X.  Guan,
``Feedback-based target localization in underwater sensor networks: A multisensor fusion approach,''
{\it  IEEE Trans.   Signal Inf.   Process.   Netw.},  vol.5,  no.1,  pp. 168--180, March 2019  doi: 10.1109/TSIPN.2018.2866335.


\bibitem{Huang2020}
S.  Huang,   P.  Chou,  X.  Jin,  Y.  Zhang,  Q.  Jiang,  S.  Yao,    
``Multi-sensor image fusion using optimized support vector machine and multiscale weighted principal component analysis,''
{\it  Electron.}, {vol. 9},  p.1531,  Aug.  2020,   doi:10.3390/electronics9091531.



\bibitem{Liping2021}
Y.  Liping, Y.  Xia, L.  Jiang,   
{\it Multisensor Fusion Estimation Theory and Application},  
1$^{st}$ ed.,  Singapore:  Springer,    2021.











%


\bibitem{Hu2020}
Z.  Hu,  J.  Hu,  G. Yang,
``A survey on distributed filtering, estimation and fusion for nonlinear systems with communication constraints: new advances and prospects,''
{\it Syst. Sci. Control Eng. }, vol.8, no.1,  pp. 189--205,  Mar.  2020,  doi:10.1080/21642583.2020.1737846.



\bibitem{Li2022}
L. Li, M. Fan, Y. Xia and Q. Geng,
``Dynamic event-triggered feedback fusion estimation for nonlinear multi-sensor systems with auto/cross-correlated Noises,''
{\it  IEEE Trans.   Signal Inf.   Process.   Netw.},  vol.8,  pp. 868--882,  Sept.  2022,  doi:10.1109/TSIPN.2022.3211172.




\bibitem{Jin2022}
H.  Jin,  S.  Sun,  
``{Distributed filtering for sensors networks with fading measurements and compensations for transmission delays and losses,}'' 
{\it Signal Process.},  {vol.  190},  p.  108306,  Jan.2022,  doi:10.1016/j.sigpro.2021.108306.



%
%
%



\bibitem{Sun2018}
S.  Sun,  F.  Peng,  H.  Lin, 
``Distributed asynchronous fusion estimator for stochastic uncertain systems with multiple sensors of different fading measurement rates, ''
{\it IEEE Trans. Signal Process. },  vol. 66, no.3,  pp.641--653,  Mar.2018,  doi:10.1109/TSP.2017.2770102. 



\bibitem{Hu2024}
J.  Hu, Z. Hu,  R.  Caballero-Águila,  X.  Yi,
``Distributed fusion filtering for multi-sensor nonlinear networked systems with multiple fading measurements via stochastic communication protocol,''
{\it Inf.  Fusion},  vol.112,  p.102543,  Dec.2024,  doi:10.1016/j.inffus.2024.102543.


\bibitem{Jin2024}
H.  Jin, Z. Du, J.  Ma, 
``Distributed filtering for networked stochastic nonlinear systems with fading measurements and random packet dropouts,'' {\it IEEE Access. }, vol.12,  pp. 110260--110272,  Mar.2024,  doi:10.1109/ACCESS.2024.3439681.



\bibitem{GarciaLigero2020}
M.J. García-Ligero, A. Hermoso-Carazo, J.  Linares-Pérez, 
``Distributed fusion estimation with sensor gain degradation and Markovian delays,''
{\it Mathematics},  vol.8,  p.1948,  Nov.2020,  doi:10.3390/math8111948.  


\bibitem{Caballero2020b}
R.  Caballero-Águila,  A.  Hermoso-Carazo,   J.  Linares-Pérez,  
``Networked fusion estimation with multiple uncertainties and time-correlated channel noise,''
{\it Inf.  Fusion},  vol.  54,  pp. 161--171,  Feb.2020,  doi:10.1016/j.inffus.2019.07.008.










\bibitem{Talebi2016}
S. Talebi, S. Kanna, D.P. Mandic,    ``A distributed quaternion Kalman filter with applications to smart grid and target tracking,''
{\it  IEEE Trans.   Signal Inf.   Process.   Netw.}, vol. 2, pp. 477--48,  Apr.2016,  doi:10.1109/TSIPN.2016.2618321.
 




\bibitem{Talebi2020}
S.P. Talebi, S. Werner,  D.P.  Mandic,  
``Quaternion-valued distributed filtering and control,''
{\it IEEE Trans.   Autom.   Control.},  {vol.65},  pp.  4246--4256,  2020,  doi:10.1109/TAC.2020.3007332.  



\bibitem{Jimenez2021}
J.D.  Jim\'enez-L\'opez, R.M. Fern\'andez-Alcal\'a, J. Navarro-Moreno, J.C.  Ruiz-Molina, 
``The distributed and centralized fusion filtering problems of tessarine signals from multi-sensor randomly delayed and missing observations under $\mathbb{T}_k$-properness conditions,''  
{\it Mathematics},  {vol.9},  p.2961,  Nov.2021,  doi:10.3390/math9222961.   




\bibitem{Fernandez2023}
R.M.  Fernández-Alcalá, J.D.  Jim\'enez-L\'opez, N.  Le Bihan, C.  Cheong Took, 
``An optimal linear fusion estimation algorithm of reduced dimension for   $\mathbb{T}$-proper
systems with multiple packet dropouts,''
{\it Mathematics}, {vol.23},  p.4047, Apr.2023, doi:10.3390/s23084047.

\bibitem{Farkas2024}
M. Farkas,  S. Rózsa,  B.  Vanek, 
``Multi-sensor attitude estimation using quaternion constrained GNSS ambiguity resolution and dynamics-based observation synchronization,''
{\it Acta Geod. Geophys. },  vol.59,  pp. 51--71,   Apr.2024,   doi:10.1007/s40328-024-00441-2.




\bibitem{Valle2020}
M.E. Valle and R.A. Lobo,  ``Quaternion-valued recurrent projection neural networks on unit quaternions, '' {\it Theor.  Comput.  Sci.}, vol.843, pp. 136--152,  Dec.2020, doi:10.1016/j.tcs.2020.08.033.


\bibitem{Valle2021}
M.E. Valle,  R.A. Lobo, 
``Hypercomplex-valued recurrent correlation neural networks,''
{\it Neurcomputing},  vol.432,  pp. 111--123,  Apr.2021,  doi:10.1109/MSP.2024.3365463.

\bibitem{Vieira2022}
G. Vieira and M.E. Valle,  ``A general framework for hypercomplex-valued extreme learning machines, '' {\it J. Comput. Math. Data Sci.},  vol. 3,  p. 100032,  Jun.2022,   doi:10.1016/j.jcmds.2022.100032.


\bibitem{Grassucci2021}
E. Grassucci, D. Comminiello, A. Unicini, 
``An information-theoretic perspective on proper quaternion variational autoencoders,''
{\it Entropy},  vol.23, no.7, p.  856,  Jul.2021,  doi:10.1109/MSP.2024.3365463.

 
\bibitem{Cariow2021}
A. Cariow, G. Cariowa,
``Fast algorithms for quaternion-valued convolutional neural networks,''
{\it  IEEE Trans. Neural Netw.  Learn. Syst. },  vol.32,  no.1,  pp.457--462,   Jan.2021,  doi:10.1109/TNNLS.2020.2979682.

\bibitem{Guizzo2023}
E.  Guizzo, T.  Weyde,  S.  Scardapane,  D.  Comminiello, 
``Learning speech emotion representations in the quaternion domain,''
{\it IEEE/ACM Trans. Audio Speech Lang.  Process.},  vol.31, pp. 1200--1212,  Mar.2023,   doi:10.1109/TASLP.2023.3250840. 



\bibitem{Clive2024}
C.  Cheong Took,  S.P.  Talebi,  R.M. Fern\'andez Alcal\'a,  D.P.  Mandic, 
``Augmented statistics of quaternion random variables: A lynchpin of quaternion learning machines,''
{\it IEEE Signal Process. Mag. }, vol.41, no.3, pp. 72--87,  May2024,  doi:10.1109/MSP.2024.3384178.

\bibitem{Cariow2024}
A. Cariow, G. Cariowa,
``Reduced-complexity algorithms for tessarine neural networks,''
{\it IEEE Trans. Neural Netw.  Learn. Syst. },  Apr.2024,  doi: 10.1109/TNNLS.2024.3385120.  Epub ahead of print. PMID: 38619958.

\bibitem{Kosal2024a}
H.H. Kosal, E. Kisi, M. Akyigit and B. Celik,  ``Elliptic quaternion matrices: theory and algorithms, '' {\it Axioms},  vol.13, no.10, p.656,   Sep.2024,  doi:10.3390/axioms13100656.

\bibitem{Kosal2024b}
H.H. Kosal, E. Kisi, M. Akyigit and B. Celik, Elliptic quaternion matrices: a MATLAB Toolbox and applications for image processing, Axioms 13~(11) (2024) 771.




%
%
%

\bibitem{Ortolani2019}
F.  Ortolani, D.  Comminiello,  M.  Scarpiniti, A.  Uncini, 
``On 4-dimensional hypercomplex algebras in adaptive signal processing,''
in: {\it Neural Advances in Processing Nonlinear Dynamic Signals},  A. Esposito, M. Faundez-Zanuy, F. Morabito, E Pasero,  Eds.  Smart Innovation, Systems and Technologies,  vol.102,  
Berlin/Heidelberg: Springer,  2019, pp. 131--140.



\bibitem{Pei2004} 
S.H. Pei,  J.H. Chang,  J.J. Ding, 
``Commutative reduced biquaternions and their Fourier transform for signal and image processing applications, ''
{\it  IEEE Trans. Signal Process. },  vol.52, no.7,  pp. 2012--2031,  Jul.2004.

\bibitem{Pei2008} 
S.H. Pei,  J.H. Chang, J.J. Ding, M.Y. Chen,  
``Eigenvalues and singular value decompositions of reduced biquaternion matrices, ''
{IEEE Trans. Circuits Syst. I Regul. Pap.},  vol.525, pp. 2673--2685,  Nov.2008, doi:10.1109/TCSI.2008.920068.

\bibitem{Melegy2022}
M.T. El-Melegy, A.T. Kamal, 
``Linear regression classification in the quaternion and
reduced biquaternion domains, ''
{\it  IEEE Signal Process Lett.},  vol.29, pp. 469--473,  Apr.2022, doi:10.1109/LSP.2022.3140682.

\bibitem{Guo2024}
Zhenwei Guo, Tongsong Jiang, Gang Wang, V.I. Vasil’ev, 
``Algebraic algorithms for eigen-problems of a reduced biquaternion matrix and applications,''
{\it  Appl. Math. Comput. }, vol. 463,  p.128358,  Feb.2024, doi:10.1016/j.amc.2023.128358.

\bibitem{Gai2024}
Z. Gai,  X. Huang,
 ``Regularization method for reduced biquaternion neural network,''
{\it Appl. Soft Comput. }, vol.166,  p.112206,  Sep.2024, doi:10.1016/j.asoc.2024.112206.


\bibitem{Nitta2019}
T.  Nitta, M.  Kobayashi,  D.P.  Mandic, 
``Hypercomplex widely linear estimation through the lens of underpinning geometry,''
{\it  IEEE Trans.   Signal Process.},  {vol.67}, pp. 3985--3994,  Jun.2019, doi:10.1109/TSP.2019.2922151.    


\bibitem{Navarro2021}
J.   Navarro-Moreno,   J.C.   Ruiz-Molina, 
``Wide-sense Markov signals on the tessarine domain.   A study under properness conditions,''
 {\it Signal Process.}, vol. 183,  p.108022,  Aug.2021, doi:10.1016/j.sigpro.2021.108022.




\bibitem{Fernandez2019}
R.M.  Fern\'andez-Alcal\'a, J. D.  Jiménez-López,  J.  Navarro-Moreno,  J.C.  Ruiz-Molina,  N. Le Bihan,
``Estimation of Widely Factorizable Hypercomplex Signals with Uncertain Observations,''
in: {\it 2019 IEEE International Conference on Acoustics, Speech and Signal Processing (ICASSP)},  Brighton, UK,  2019, pp. 8504-8508. 


\bibitem{Took18}
M.  Xiang; S.  Enshaeifar,  A.E.  Stott, C.  Cheong Took,   Y.  Xia, S.  Kanna,  D.P.  Mandic,   
``Simultaneous diagonalisation of the covariance and complementary covariance matrices in quaternion widely linear signal processing,''
{\it Signal Process.},  {vol.148}, pp. 193--204,  {Jan.2018},   doi:10.48550/arXiv.1705.00058.

\bibitem{Navarro2022}
J.  Navarro-Moreno,  R.M.  Fernández-Alcalá,  J.C.  Ruiz-Molina,  
``Proper ARMA modeling and forecasting in the generalized Segre’s quaternions domain,''
{\it Mathematics},  vol.10,  p.1083,  Mar.2022,  doi:10.3390/math10071083.


\bibitem{Navarro2024}
J. D.  Jiménez-López, J.  Navarro-Moreno, R. M.  Fernández-Alcalá,  J. C.  Ruiz-Molina, 
``Proper processing of  $\beta$-quaternion wide-sense Markov signals from randomly lost observations,''
{\it IEEE Signal Process  Lett. },  vol.31, pp. 2060--2064,   Jan.2024,  doi:10.1109/LSP.2024.3440964. 


\end{thebibliography}


\begin{IEEEbiographynophoto}
{Rosa M. Fernández-Alcalá} (Jaén, Spain)  received the M.Sc. degree in mathematics from the University of Granada, Spain, in 1996, and the Ph.D. degree in mathematics from the University of Jaén,  Spain, in 2002.

She has been with the Department of Statistics and Operations Research, University of Jaén, since 1996, where she is currently a Professor. Her research interests include widely linear estimation methods and proper hypercomplex signal processing as a tool for dimensional reduction. She has published numerous research papers and has participated in competitive research projects under the (I+D+i) National Research Program, co-funded by the European Regional Development Funds (ERDF).

Prof.  Fernández-Alcalá is a member of the research group “Estadística Teórica y Aplicada e Investigación Operativa” at the University of Jaén. 

\end{IEEEbiographynophoto}

\begin{IEEEbiographynophoto}
{José D. Jiménez-López } was born in Granada, Spain.    He received the M.Sc. degree in mathematics and the Ph.D. degree from the University of Granada, Spain, in 1997 and 2000, respectively.

In 1999,  he joined the Department of Statistics and Operations Research, University of Jaén,  Spain, where he has been developing research and teaching activities,  and since 2009,  as Associate Professor.  He has participated as researcher in different competitive National Program Research Projects.  His research has focused on elliptical contoured multivariate distributions,  and he is currently working on hypercomplex signal estimation problems.

Dr. Jiménez-López is associated with the research group “Estadística Teórica y Aplicada e Investigación Operativa” at the University of Jaén.
\end{IEEEbiographynophoto}

\begin{IEEEbiographynophoto}
{Jesús Navarro-Moreno} was born in Jaén,  Spain, in 1970. He received the M.Sc. degree in mathematics and the Ph.D. degree from the University of Granada, Spain, in 1993 and 1998, respectively.

In 1993, he joined the faculty at the University of Jaén, where he has been a Professor in the Department of Statistics and Operations Research since October 2018. His research interests include second-order stochastic processes, statistical signal processing, and estimation and detection theory.

Prof.   Navarro-Moreno is a member of the Spanish Society of Statistics and Operations Research (SEIO).

\end{IEEEbiographynophoto}

\begin{IEEEbiographynophoto}
{Juan Carlos Ruiz-Molina } was born in Jaén, Spain, in 1966. He received the M.Sc. degree in mathematics and the Ph.D. degree from the University of Granada, Spain, in 1989 and 1993, respectively.

In 1989, he joined the faculty at the University of Granada, where he held various teaching and research positions in the Department of Statistics and Operations Research from 1989 to 1993. Since 2003, he has been a Professor in the Department of Statistics and Operations Research at the University of Jaén, Spain. His research interests include identification, estimation, and detection. He has also authored a textbook on estimation in Spanish.

Prof.  Ruiz-Molina is part of the research group “Estadística Teórica y Aplicada e Investigación Operativa” at the University of Jaén. 
\end{IEEEbiographynophoto}

%

\end{document}